\documentclass[11pt]{article}
\usepackage[english]{babel}
\usepackage{amsmath,amsthm,amssymb}
\usepackage{mathrsfs}
\usepackage{bbm}
\usepackage{amsfonts}
\usepackage[margin=0.8in]{geometry}
\usepackage{graphicx}
\usepackage{caption}
\usepackage{subcaption}
\usepackage{wasysym}
\usepackage{enumerate}
\usepackage{algorithm2e}
\usepackage{tabularx}
\usepackage{color}
\usepackage{wrapfig}
\usepackage{todonotes}

\newtheorem{thm}{Theorem}[section]
\newtheorem{ass}[thm]{Assumption}
\newtheorem{cor}[thm]{Corollary}
\newtheorem{lem}[thm]{Lemma}
\newtheorem{prop}[thm]{Proposition}
\newtheorem*{hyp*}{Hypothesis}
\theoremstyle{definition}
\newtheorem{defn}[thm]{Definition}

\theoremstyle{rem}
\newtheorem{rem}[thm]{Remark}

\numberwithin{equation}{section}

\newcommand{\R}{\mathbb R}
\newcommand{\eps}{\varepsilon}

\newcommand{\bbF}{\mathbb F}

\newcommand{\bbN}{\mathbb N}

\newcommand{\mcA}{\mathcal{A}}

\newcommand{\mcB}{\mathcal{B}}

\newcommand{\mcL}{\mathcal L}
\newcommand{\mcE}{\mathcal E}

\newcommand{\mcF}{\mathcal F}

\newcommand{\mcV}{\mathcal V}
\newcommand{\mcT}{\mathcal T}

\newcommand{\mcM}{\mathcal M}
\newcommand{\mcH}{\mathcal H}
\newcommand{\mcK}{\mathcal K}
\newcommand{\mcU}{\mathcal U}
\newcommand{\mcS}{\mathcal S}

\newcommand{\Pred}{\mathcal{P}}

\newcommand{\E}{\mathbb{E}}
\newcommand{\Prob}{\mathbb{P}}

\DeclareMathOperator*{\esssup}{ess\,sup}
\newcommand{\essinf}{\mathop{\rm{ess}\inf}}

\newcommand{\ett}{\mathbbm{1}}

\newcommand{\cadlag}{c\`adl\`ag~}

\newcommand{\trace}{{\rm Tr}}

\newcommand{\ie}{\textit{i.e.\ }}
\newcommand{\eg}{\textit{e.g.\ }}

\newcommand{\rev}{}

\begin{document}

\title{Optimal Stopping of BSDEs with Constrained Jumps and Related Double Obstacle PDEs}

\author{Magnus Perninge\footnote{M.\ Perninge is with the Department of Physics and Electrical Engineering, Linnaeus University, V\"axj\"o,
Sweden. e-mail: magnus.perninge@lnu.se.}} %
\maketitle
\begin{abstract}
We consider partial differential equations (PDEs) characterized by an upper barrier that depends on the solution itself and a fixed lower barrier, while accommodating a non-local driver. First, we show a Feynman-Kac representation for the PDE when the driver is local. Specifically, we relate the non-linear Snell envelope, arising from an optimal stopping problem--where the underlying process is the first component in the solution to a stopped backward stochastic differential equation (BSDE) with jumps and a constraint on the jumps process--to a viscosity solution for the PDE. Leveraging this Feynman-Kac representation, we subsequently prove existence and uniqueness of viscosity solutions in the non-local setting by employing a contraction argument. This approach also introduces a novel form of non-linear Snell envelope and expands the probabilistic representation theory for PDEs.
\end{abstract}

\section{Introduction}
We consider partial differential equations (PDEs) of the type
\begin{align}\label{ekv:var-ineq}
\begin{cases}
  \min\{v(t,x)-h(t,x),\max\{v(t,x)-\mcM v(t,x),-v_t(t,x)-\mcL v(t,x)\\
  \quad-f(t,x,v(t,\cdot),\sigma^\top(t,x)\nabla_x v(t,x))\}\}=0,\quad\forall (t,x)\in[0,T)\times \R^d \\
  v(T,x)=\psi(x),
\end{cases}
\end{align}
where for each $(t,x,z)\in [0,T]\times\R^d\times\R^d$, the map $g\mapsto f(t,x,g,z):C(\R^d\to\R)\to\R$ is a functional, $\mcM v(t,x):=\inf_{e\in E}\{v(t,x+\gamma(t,x,e))+\chi(t,x,e)\}$ and
\begin{align}
  \mcL:=\sum_{j=1}^d a_j(t,x)\frac{\partial}{\partial x_j}+\frac{1}{2}\sum_{i,j=1}^d (\sigma\sigma^\top(t,x))_{i,j}\frac{\partial^2}{\partial x_i\partial x_j}
\end{align}
is the infinitesimal generator related to the stochastic differential equation (SDE)
\begin{align*}
\check X_s=x + \int_0^s a(r,\check X_r)dr + \int_0^s \sigma(r,\check X_r)dW_r,\quad\forall s\in[0,T].
\end{align*}

We establish existence and uniqueness (within the set of continuous functions of polynomial growth) of viscosity solutions to the above PDE by relating $v(t,x)$ to $Y^{t,x}_t$, where for each $(t,x)\in [0,T]\times\R^d$,
\begin{align}\label{ekv:opt-stop}
  Y^{t,x}_s=\esssup_{\tau\in\mcT_s}Y^{t,x,\tau}_s,\quad\forall s\in [t,T]
\end{align}
and for each stopping time $\tau\geq t$, the process the $Y^{t,x,\tau}$ is the first component in the quadruple of processes $(Y^{t,x,\tau},Z^{t,x,\tau},V^{t,x,\tau},K^{t,x,\tau})$ that is the unique maximal solution to the backward stochastic differential equation (BSDE), driven by the above mentioned Brownian motion and an independent Poisson random measure $\mu$ with intensity $\lambda$, that has a constraint on the jump component,
\begin{align}\label{ekv:stopped-bsde}
  \begin{cases}
     Y^{t,x,\tau}_s=\Psi(\tau,X^{t,x}_\tau)+\int_s^\tau f(r,X^{t,x}_r, \bar Y(r,\cdot), Z^{t,x,\tau}_r)dr -\int_s^\tau  Z^{t,x,\tau}_r dW_r-\int_{s}^\tau\!\!\!\int_E  V^{t,x,\tau}_{r}(e)\mu(dr,de)
    \\
    \quad-(K^{-,t,x,\tau}_\tau- K^{-,t,x,\tau}_s),\quad \forall s\in [t,\tau],
    \\
    V^{t,x,\tau}_s(e)\geq - \chi(s,X^{t,x}_{s-},e),\quad d\Prob\otimes ds\otimes \lambda(de)-a.e.
  \end{cases}
\end{align}
Here, $\Psi(t,x):=\ett_{[t<T]}h(t,x)+\ett_{[t=T]}\psi(x)$,  the process $X^{t,x}$ is the unique solution to the forward SDE
\begin{align}\label{ekv:SDE-w-jmp}
X^{t,x}_s=x + \int_t^s a(r,X^{t,x}_r)dr + \int_t^s \sigma(r,X^{t,x}_r)dW_r+\int_{t}^s\!\!\!\int_E  \gamma(r,X^{t,x}_{r-},e)\mu(dr,de),\quad\forall s\in[0,T]
\end{align}
and the driver at time $t\in[0,T]$ depends on the solution to the optimal stopping problem through the map $(t,x)\mapsto \bar Y(t,x):=Y^{t,x}_t$. Our contribution is twofold in the sense that our approach implicitly provides existence of a unique solution to the optimal stopping problem \eqref{ekv:opt-stop}-\eqref{ekv:SDE-w-jmp}.\\

\textbf{Related literature} Partial differential equations featuring non-local drivers, wherein the mapping $g \mapsto f(t,x,g,z)$ adopts the structure of an integral, are commonly referred to as integro-partial differential equations (IPDEs). A branch of the literature that relates BSDEs with jumps to IPDEs have focused on the setting where the driver takes the form
\begin{align*}
f(t,x,g,z)=\bar f(t,x,g(x),z,\int_E\theta(t,x,e)(g(x+\zeta(t,x,e))-g(x))\lambda(de)),
\end{align*}
for a Lipschitz function $\bar f$. In particular, \cite{Barles1997} and \cite{Dumitrescu15} both assume that $\bar f$ is non-decreasing in the last variable and that $\theta\geq 0$, the former considering the relation between regular BSDEs with jumps and IPDEs, whereas the latter derive a relation between the reflected BSDEs with jumps considered in \cite{Quenez2014} and IPDEs with one-sided obstacles. These results where later extended in \cite{HamMor16} to the setting when $\bar f$ is a general Lipschitz function and $\theta$ may be negative. We also mention the work in \cite{HamMnif2023}, where a system of IPDEs with interconnected obstacles are treated.

When the lower barrier $h$ is absent, the PDE described in \eqref{ekv:var-ineq} is termed a quasi-variational inequality (QVI). It is well known that, under suitable conditions on the involved parameters, value functions to impulse control problems are solutions (in viscosity sense) to standard QVIs when the driving noise process is a Brownian motion (see the seminal work in~\cite{BensLionsImpulse}) and to so called quasi-integrovariational inequalities when the driving noise is a general L\'evy process~\cite{OksenSulemBok}. Employing a contraction argument akin to the one delineated in the present study, QVIs featuring general non-local drivers were investigated in \cite{qvi-rbsde}. Specifically, \cite{qvi-rbsde} establishes a Feynman-Kac representation for such QVIs by linking their solutions to systems of reflected BSDEs (RBSDEs).

An alternative Feynman-Kac representation for solutions to standard QVIs was proposed in \cite{KharroubiQVI} (see also \cite{Bouchard09}), where the solution to a QVI is related to the minimal solution of a BSDE with constrained jumps. In particular, their result implies that (when $f$ is local) the deterministic function $v$ defined through the relation $v(t,x):=Y^{t,x,T}_t$, with $Y^{t,x,T}$ as in \eqref{ekv:stopped-bsde}, is the unique viscosity solution to the PDE obtained by setting $h\equiv-\infty$ in \eqref{ekv:var-ineq}.

The ensemble of approaches to find probabilistic representations of PDEs or solve stochastic optimal control problems that utilize BSDEs with constrained jumps is commonly referred to as \emph{control randomization}. A significant breakthrough in this field was achieved in \cite{Fuhrman15}, which directly linked the value function of the randomized control problem to that of the original control problem. This eliminated the need for a Feynman-Kac representation, thereby expanding the theoretical framework to encompass stochastic systems with path-dependencies. Building upon this foundation, subsequent advancements extended their approach to the framework of partial information settings in \cite{Bandini18} and optimal switching problems in \cite{Fuhrman2020}. Recently, \cite{imp-stop-game} further extended the scope of control randomization to zero-sum games by (within a non-markovian framework) relating the solution to the above optimal stopping problem \eqref{ekv:opt-stop}-\eqref{ekv:SDE-w-jmp} to the upper and lower value functions in a stochastic differential game between an impulse controller and a stopper. Notable is that \cite{imp-stop-game} presumes that $f$ only depends on local values of $Y$, whereas it may depend on $V$.

An intermediate result in the present work addresses the framework of a local driver and bridges a gap left in \cite{imp-stop-game}. Specifically, we establish a connection between the non-linear Snell envelope examined in \cite{imp-stop-game} and viscosity solutions to \eqref{ekv:var-ineq}. Consequently, we elucidate that, in the Markovian framework, the value function of the aforementioned zero-sum game indeed constitutes the unique viscosity solution to \eqref{ekv:var-ineq}. Moreover, our primary findings extend those of \cite{imp-stop-game} in the Markovian framework by proving the existence of the more general non-linear Snell envelope defined through equations \eqref{ekv:opt-stop}-\eqref{ekv:SDE-w-jmp}.\\

\textbf{Outline} In the next section, we set the notations and state the assumptions that hold throughout. In addition, we give some preliminary results that are repeatedly referred in the article. Then, in Section~\ref{sec:local} we turn to the local setting before we, in the following section, derive the complete result. Uniqueness of solutions to~\eqref{ekv:var-ineq}, when the driver $f$ is local, appears rudimentary and resembles earlier results deduced in, for example, \cite{Morlais13}. However, since our setting is fundamentally different and for the sake of completeness, a uniqueness proof through viscosity comparison is included as an appendix.

\section{Preliminaries\label{sec:prel}}
\subsection{Notation}
We let $(\Omega,\mcF,\Prob)$ be a complete probability space on which lives a $d$-dimensional Brownian motion $W$ and an independent Poisson random measure $\mu$ defined on $[0,T]\times E$ with intensity $\lambda(de)$. Here, it is assumed that $E$ is a compact subset of $\R^d$ endowed with its Borel $\sigma$-field $\mcB(E)$. We denote by $\bbF:=(\mcF_t)_{0\leq t\leq T}$ the augmented natural filtration generated by $W$ and $\mu$ and for $t\in[0,T]$ we let $\bbF^t:=(\mcF^t_s)_{t\leq s\leq T}$ (resp. $\bbF^{W,t}:=(\mcF^{W,t}_s)_{t\leq s\leq T}$) denote the augmented natural filtration generated by $(W_s-W_t:t\leq s\leq T)$ and $\mu(\cdot\cap[t,T],\cdot)$ (resp. $(W_s-W_t,t\leq s\leq T)$).\\

\noindent Throughout, we will use the following notation, where $d\geq 1$ is the dimension of the state-space:
\begin{itemize}
\item We let $\Pi^g$ denote the set of all functions $\varphi:[0,T]\times\R^d\to\R$ that are of polynomial growth in $x$, \ie there are constants $C,\rho>0$ such that $|\varphi(t,x)|\leq C(1+|x|^\rho)$ for all $(t,x)\in [0,T]\times\R^d$, and let $\Pi^g_c$ be the subset of jointly continuous functions.
  \item For each $t\in[0,T]$, $\Pred_t$ is the $\sigma$-algebra of $\bbF^t$-predictably measurable subsets of $[t,T]\times \Omega$ and $\Pred:=\Pred_0$.
  \item We let $\mcT$ be the set of all $[0,T]$-valued $\bbF$-stopping times and for each $\eta\in\mcT$, we let $\mcT_\eta$ be the subset of stopping times $\tau$ such that $\tau\geq \eta$, $\Prob$-a.s. Moreover, we let $\mcT^t$ (resp. $\mcT^t_\eta$) be the corresponding subsets of $\bbF^t$-stopping times, with $\tau\geq t$ (resp. $\tau\geq \eta$), $\Prob$-a.s.
  \item For $p\geq 1$, $t\in [0,T]$ and $\tau\in\mcT^t$, we let $\mcS^{p}_{[t,\tau]}$ be the set of all $\R$-valued, $\bbF^t$-progressively measurable, \cadlag processes $(Z_s: s\in [t,\tau])$ for which $\|Z\|_{\mcS^p_{[t,\tau]}}:=\E\big[\sup_{s\in[t,\tau]} |Z_s|^p\big]<\infty$. Moreover, we let $\mcS^p_{[t,\tau],i}$ be the subset of $\bbF^t$-predictably measurable and non-decreasing processes with $Z_t=0$. Whenever $\tau=T$ we use the notations $\mcS^p_{t}$ and $\mcS^p_{t,i}$, respectively.
  \item We let $\mcH^{p}_{[t,\tau]}(W)$ denote the set of all $\R^d$-valued $\bbF^t$-progressively measurable processes $(Z_s: s\in[t,\tau])$ such that $\|Z\|_{\mcH^p_{[t,\tau]}(W)}:=\E\big[\big(\int_t^\tau |Z_s|^2 ds\big)^{p/2}\big]^{1/p}<\infty$. Furthermore, we set $\mcH^{p}_{t}(W):=\mcH^{p}_{[t,T]}(W)$.
  \item We let $\mcH^{p}_{[t,\tau]}(\mu)$ denote the set of all $\R$-valued, $\Pred_t \otimes\mcB(E)$-measurable maps $(Z_s(e): s\in[t,\tau],e\in E)$ such that $\|Z\|_{\mcH^p_{[t,\tau]}(\mu)}:=\E\big[\big(\int_t^\tau \!\!\int_E |Z_s(e)|^2 \lambda(de)ds\big)^{p/2}\big]^{1/p}<\infty$ and set $\mcH^{p}_{t}(\mu):=\mcH^{p}_{[t,T]}(\mu)$.
  \item For $t\in[0,T]$, we let $\mcA_t$ denote the set of all $[-1,1]^d$-valued, $\bbF^t$-progressively measurable processes $(\alpha_s:t\leq s\leq T)$ and set $\mcA:=\mcA_0$. Moreover, we let $\mcA^W_t$ be the subset of $\bbF^{W,t}$-progressively measurable processes (resp. $\mcA^W:=\mcA^W_0$).
  \item For $t\in[0,T]$, we define the composition $\oplus_t$ of $\alpha^1\in\mcA$ and $\alpha^2\in\mcA_t$ as $(\alpha^1\oplus_t\alpha^2)_s:=\ett_{[0,t)}(s)\alpha^1_s+\ett_{[t,T]}(s)\alpha^2_s$.
  \item We let $\mcV_t$ denote the set of all $\Pred_t\otimes\mcB(E)$-measurable, bounded maps $\nu:[t,T]\times\Omega\times E\to [0,\infty)$ and for each $n\in\bbN$, we denote by $\mcV^n_t$ the subset of maps $\nu:[t,T]\times\Omega\times E\to [0,n]$.
\end{itemize}

We also mention that, unless otherwise specified, all inequalities between random variables are to be interpreted in the $\Prob$-a.s.~sense.

\subsection{Assumptions}
We make the following assumption on the intensity $\lambda$ of the process $\mu$.
\begin{ass}\label{ass:on-lambda}
We assume that $\lambda$ has full topological support on $E$ and finite intensity, \ie $\lambda(E)<\infty$.
\end{ass}
Throughout, we make the following assumptions on the parameters, where $\rho>0$ is a fixed constant:
\begin{ass}\label{ass:on-coeff}
\begin{enumerate}[i)]
  \item\label{ass:on-coeff-f} We assume that $f:[0,T]\times \R^d\times C(\R^d\to\R)\times\R^{d}\to \R$ ($(t,x,g,z)\to f(t,x,g,z)$) is such that for any $v\in\Pi^g_c$, the map $(t,x)\mapsto f(t,x,v(t,\cdot),z)$ is jointly continuous, uniformly in $z$, $f$ is of polynomial growth in $x$, \ie there is a $C_f>0$ such that
  \begin{align*}
    |f(t,x,0,0)|\leq C_f(1+|x|^\rho)
  \end{align*}
  and that there are constants $k_f,K_\Gamma>0$ such that for any $t\in[0,T]$, $x\in\R^d$, $g,\tilde g\in C(\R^d\to\R)$ and $z,\tilde z\in\R^{d}$ we have
  \begin{align}\label{ekv:f-lipschitz}
    |f(t,x,\tilde g,\tilde z)-f(t,x,g,z)|&\leq k_f(\sup_{x'\in \Lambda_f(|x|)}|\tilde g(x')-g(x')|+|\tilde z-z|),
  \end{align}
  where for each $\alpha\in\R_+$, $\Lambda_f(\alpha):=\{x\in\R^d:\|x\|\leq \alpha\vee K_\Gamma\}$ is the closed ball of radius $\alpha\vee K_\Gamma$ centered at the origin.
  \item\label{ass:on-coeff-h} The lower barrier $h:[0,T]\times \R^d \to \R$ is jointly continuous and of polynomial growth in $x$, \ie there is a $C_h>0$ such that
  \begin{align*}
    |h(t,x)|\leq C_h(1+|x|^\rho).
  \end{align*}
  \item\label{ass:on-coeff-psi} The terminal reward $\psi:\R^d\to\R$ is continuous and satisfies the growth condition
  \begin{align*}
    |\psi(x)|\leq C_\psi(1+|x|^\rho)
  \end{align*}
  for some $C_\psi>0$.
  \item\label{ass:on-coeff-chi} The jump-barrier $\chi:[0,T]\times \R^d\times E\to \R_+$ is jointly continuous and of polynomial growth, \ie
   \begin{align*}
    |\chi(t,x,e)|\leq C_\chi(1+|x|^\rho)
  \end{align*}
  for some $C_\chi>0$.
  \item\label{ass:no-free-loop} \rev{\underline{No-free-$\eps$-loop:} There are functions $\delta_1,\delta_2:[0,T]\to (0,\infty)$ such that whenever there is a $x_0\in\R^d$, a $k\in\bbN$ and a sequence $(e_j)_{j=0}^{k-1}$ in $E$ such that for some $t\in[0,T]$, we have $\|x_k-x_0\|\leq \delta_1(t)$ where $x_j=x_{j-1}+\gamma(t,x_{j-1},e_{j-1})$ for $j=1,\ldots,k$, then
  \begin{align*}
    \sum_{j=0}^{k-1}\chi(t,x_{j},e_j)\geq\delta_2(t).
  \end{align*}}
  \item\label{ass:on-coeff-@end} For each $(x,e)\in\R^d\times E$, we have
  \begin{align*}
    h(T,x)\leq \psi(x)\leq \psi(x+\gamma(T,x,e))+\chi(T,x,e).
  \end{align*}
\end{enumerate}
\end{ass}

\rev{Except for \eqref{ekv:f-lipschitz},} the above assumptions are all standard in the context of impulse control and for the related QVIs. \rev{The general type of Lipschitz continuity imposed on the map $g\mapsto f(t,x,g,z):C(\R^d\to\R)\to\R$ by \eqref{ekv:f-lipschitz} allows for modeling scenarios where, for instance, $f(t,x,v(t,\cdot),z)=\tilde f(t,x,v(t,x),z)+k\mcM v(t,x)$ for some constant $k\in\R$. From an application perspective, this corresponds to impulse control problems where the running cost is influenced by the optimal remaining cost that one could obtain by invoking an impulse.}

\rev{Notable is that, rather than setting a positive lower bound for the intervention cost $\chi$, the ``no-free-$\eps$-loop''-condition (\ref{ass:no-free-loop}) dictates that there is a positive cost attached to making a series of interventions that almost take the state back to the initial one. This condition generalizes the commonly applied \emph{no-free-loop} condition of optimal switching (see \eg \cite{Morlais13}).}\\

Moreover, we make the following assumptions on the coefficients of the forward SDE:

\begin{ass}\label{ass:onSDE}
For any $t,t'\in [0,T]$, $e\in E$ and $x,x'\in\R^d$ we have:
\begin{enumerate}[i)]
  \item\label{ass:onSDE-Gamma} The function $\gamma:[0,T]\times\R^d\times E\to\R^d$ is jointly continuous and satisfies the growth condition
  \begin{align}\label{ekv:imp-bound}
    |x+\gamma(t,x,e)|\leq K_\Gamma\vee |x|.
  \end{align}
  Moreover, it is Lipschitz continuous in $x$ uniformly in $(t,e)$, \ie
  \begin{align*}
    |\gamma(t,x',e)-\gamma(t,x,e)|\leq k_\gamma|x'-x|
  \end{align*}
  for some $k_\gamma>0$ and all $(t,x,x',e)\in [0,T]\times \R^d\times\R^d\times E$.
  \item\label{ass:onSDE-a-sigma} The coefficients $a:[0,T]\times\R^d\to\R^{d}$ and $\sigma:[0,T]\times\R^d\to\R^{d\times d}$ are jointly continuous and satisfy the growth conditions
  \begin{align*}
    |a(t,x)|+|\sigma(t,x)|&\leq C_{a,\sigma}(1+|x|),
  \end{align*}
  for some $C_{a,\sigma}>0$ and the Lipschitz continuity
  \begin{align*}
    |a(t,x)-a(t,x')|+|\sigma(t,x)-\sigma(t,x')|&\leq k_{a,\sigma}|x'-x|,
  \end{align*}
  for some $k_{a,\sigma}>0$.
\end{enumerate}
\end{ass}

\rev{The growth condition in \eqref{ekv:imp-bound} gives the constant $K_\Gamma$ the role of a barrier such that impulses affect the state in a non-expansive manner whenever the magnitude of the state exceeds $K_\Gamma$. This guarantees that, regardless of the choice of $u\in\mcU$, the state process $X^{t,x,u}$ does not blow up in a finite time.}

\subsection{Viscosity solutions}

We define the upper, $v^*$, and lower, $v_*$ semi-continuous envelope of a function $v:[0,T]\times\R^d\to\R$ as
\begin{align*}
v^*(t,x):=\limsup_{(t',x')\to(t,x),\,t'<T}v(t',x')\quad {\rm and}\quad v_*(t,x):=\liminf_{(t',x')\to(t,x),\,t'<T}v(t',x').
\end{align*}
Next, we introduce the limiting parabolic superjet $\bar J^+v$ and subjet $\bar J^-v$.
\begin{defn}\label{def:jets}
Subjets and superjets
\begin{enumerate}[i)]
\item For a l.s.c. (resp. u.s.c.) function $v : [0, T]\times \R^d \to \R$, the parabolic subjet, denote by $J^-v(t, x)$, (resp.~the parabolic superjet, $J^+v(t, x)$) of $v$ at $(t, x) \in [0, T]\times \R^d$, is defined as the set of triples $(p, q,M) \in \R \times\R^d \times \mathbb S^d$ satisfying
    \begin{align*}
      v(t', x') \geq (\text{resp.~}\leq)\: v(t, x) + p(t'- t)+ < q, x' - x > + \tfrac{1}{2} < x' - x,M(x' - x) > +o(|t' - t| + |x' - x|^2)
    \end{align*}
    for all $(t',x')\in[0,T]\times\R^d$, where $\mathbb S^d$ is the set of symmetric real matrices of dimension $d\times d$.
\item For a l.s.c. (resp. u.s.c.) function $v : [0, T]\times \R^d \to \R$ we denote by $\bar J^-v(t, x)$ the parabolic limiting
subjet (resp. $\bar J^+v(t, x)$ the parabolic limiting superjet) of $v$ at $(t, x) \in [0, T]\times \R^d$, defined as the set of triples $(p, q,M) \in \R \times\R^d \times \mathbb S^d$ such that:
\begin{align*}
  (p, q,M) = \lim_{n\to\infty} (p_n, q_n,M_n),\quad (t, x) = \lim_{n\to\infty}(t_n, x_n)
\end{align*}
for some sequence $(t_n,x_n,p_n,q_n,M_n)_{n\geq 1}$ with $(p_n, q_n,M_n) \in J^-v(t_n, x_n)$ (resp. $(p_n, q_n,M_n) \in J^+v(t_n, x_n)$) for all $n\geq 1$ and $v(t, x) = \lim_{n\to\infty}v(t_n, x_n)$.
\end{enumerate}
\end{defn}

We now give the definition of a viscosity solution for the QVI in \eqref{ekv:var-ineq}. (see also pp. 9-10 of \cite{UsersGuide}).
\begin{defn}\label{def:visc-sol-jets}
Let $v$ be a locally bounded function from $[0,T]\times \R^d$ to $\R$. Then,
\begin{enumerate}[a)]
  \item It is referred to as a viscosity supersolution (resp. subsolution) to \eqref{ekv:var-ineq} if it is l.s.c.~(resp. u.s.c.) and satisfies:
  \begin{enumerate}[i)]
    \item $v(T,x)\geq \psi(x)$ (resp. $v(T,x)\leq \psi(x)$)
    \item For any $(t,x)\in [0,T)\times\R^d$ and $(p,q,X)\in \bar J^- v(t,x)$ (resp. $\bar J^+ v(t,x)$) we have
    \begin{align*}
      \min\big\{&v(t,x)-h(t,x),\max\{v(t,x)-\mcM v(t,x),-p-q^\top a(t,x)
      \\
      &-\frac{1}{2}\trace(\sigma\sigma^\top(t,x)X)-f(t,x,v(t,\cdot),\sigma^\top(t,x)q)\}\big\}\geq 0\quad (\text{resp. }\leq 0)
    \end{align*}
  \end{enumerate}
  \item It is called a viscosity solution to \eqref{ekv:var-ineq} if $v_*$ is a supersolution and $v^*$ is a subsolution.
\end{enumerate}
\end{defn}

We will sometimes use the following alternative definition of viscosity supersolutions (resp. subsolutions):
\begin{defn}\label{def:visc-sol-dom}
  A l.s.c.~(resp. u.s.c.) function $v$ is a viscosity supersolution (resp.~subsolution) to \eqref{ekv:var-ineq} if $v(T,x)\leq \psi(x)$ (resp. $\geq \psi(x)$) and whenever $\varphi\in C^{1,2}([0,T]\times\R^d\to\R)$ is such that $\varphi(t,x)=v(t,x)$ and $\varphi-v$ has a local maximum (resp. minimum) at $(t,x)$, then
  \begin{align*}
    \min\big\{&v(t,x)-h(t,x),\max\{v(t,x)-\mcM v(t,x),-\varphi_t(t,x)-\mcL\varphi(t,x)
    \\
    &-f(t,x,v(t,\cdot),\sigma^\top(t,x)\nabla_x\varphi(t,x))\}\big\}\geq 0\quad(\text{resp. }\leq 0).
  \end{align*}
\end{defn}

It can readily be shown (see \eg \cite{UsersGuide}) that Definition~\ref{def:visc-sol-jets} and Definition~\ref{def:visc-sol-dom} are equivalent.


\subsection{Reflected BSDEs with jumps and obstacle problems}
We introduce the local driver $\tilde f$ that satisfies the following assumption:
\begin{ass}\label{ass:tilde-coeff}
\begin{enumerate}[i)]
  \item The driver $\tilde f:[0,T]\times\R^d\times\R\times\R^d\to\R$ is jointly continuous and of polynomial growth in $x$, \ie $|\tilde f(t,x,0,0)|\leq C(1+|x|^\rho)$. Moreover, $(t,x)\mapsto \tilde f(t,x,y,z)$ is continuous, uniformly in $(y,z)$, and $(y,z)\mapsto \tilde f(t,x,y,z)$ is Lipschitz continuous, uniformly in $(t,x)$.
\end{enumerate}
\end{ass}

We give the following proposition, strategically formulated to streamline subsequent implementation processes. It is worth noting that comparable findings are documented in \cite{Dumitrescu15,HamMor16}, with the latter being the most closely aligned with our own.

\begin{prop}\label{prop:rbsde-jmp}
For each $(t,x)\in [0,T]\times\R^d$ and $n\in\bbN$, there is a unique quadruple $(Y^{t,x},Z^{t,x},V^{t,x},K^{+,t,x})\in\mcS^{2}_t \times \mcH^{2}_t(W) \times \mcH^{2}_t(\mu) \times \mcS^{2}_{t,i}$ such that
\begin{align}\label{ekv:rbsde-jmp}
  \begin{cases}
     Y^{t,x}_s=\psi(X^{t,x}_T)+\int_s^T \tilde f^n(r,X^{t,x}_r, Y^{t,x}_r, Z^{t,x}_r, V^{t,x}_{r})dr -\int_s^T  Z^{t,x}_r dW_r-\int_{s}^T\!\!\!\int_E  V^{t,x}_{r}(e)\mu(dr,de)
    \\
    \quad+(K^{+,t,x}_T- K^{+,t,x}_s),\quad \forall s\in [t,T],
    \\
    Y^{t,x}_s\geq h(s,X^{t,x}_s),\, \forall s\in [t,T] \quad\text{and}\quad\int_t^T\!\! \big(Y^{t,x}_s-h(s,X^{t,x}_s)\big)dK^{+,t,x}_s,
  \end{cases}
\end{align}
where
\begin{align*}
\tilde f^n(t,x,y,z,v):=\tilde f(t,x,y,z)-n\int_E(v(e)+\chi(t,x,e))^-\lambda(de).
\end{align*}
Moreover, $Y^{t,x}_s=\esssup_{\tau\in\mcT_s}Y^{t,x;\tau}_s$, where for each $\tau\in\mcT_t$, the triple
$(Y^{t,x;\tau},Z^{t,x;\tau},V^{t,x;\tau})\in\mcS^{2}_{[0,\tau]} \times \mcH^{2}_{[0,\tau]}(W) \times \mcH^{2}_{[0,\tau]}(\mu)$ satisfies
\begin{align}\label{ekv:bsde-jmp}
  Y^{t,x;\tau}_s&=\Psi(\tau,X^{t,x}_\tau)+\int_s^\tau \tilde f^n(r,X^{t,x}_r, Y^{t,x;\tau}_r, Z^{t,x;\tau}_r, V^{t,x;\tau}_{r})dr
   \\
   &\quad-\int_s^\tau  Z^{t,x;\tau}_r dW_r-\int_{s}^\tau\!\!\!\int_E  V^{t,x;\tau}_{r}(e)\mu(dr,de),\quad \forall s\in [0,\tau]
\end{align}
and for each $\eta\in \mcT_t$, the stopping time
\begin{align*}
  \tau^*:=\inf\{s\geq \eta:Y^{t,x}_s=h(s,X^{t,x}_s)\}\wedge T
\end{align*}
is optimal in the sense that $Y^{t,x}_\eta=Y^{t,x;\tau^*}_\eta$ (note that $\tau^*\in\mcT^t_\eta$ whenever $\eta\in\mcT^t$). Finally, there is a function $v_n\in\Pi^g_c$ such that $v_n(s,X^{t,x}_s)=Y^{t,x}_s$ for all $s\in[t,T]$ and $v_n$ is the unique viscosity solution in $\Pi^g_c$ to
\begin{align}\label{ekv:obst-prob-n}
\begin{cases}
  \min\{v_n(t,x)-h(t,x),-\frac{\partial}{\partial t}v_n(t,x)-\mcL v_n(t,x)+\mcK^n v_n(t,x)\\
  \quad-\tilde f(t,x,v_n(t,x),\sigma^\top(t,x)\nabla_x v_n(t,x))\}=0,\quad\forall (t,x)\in[0,T)\times \R^d \\
  v_n(T,x)=\psi(x),
\end{cases}
\end{align}
where $\mcK^n \phi(t,x):=n\int_E(\phi(t,x+\gamma(t,x,e)) + \chi(t,x,e)-\phi(t,x))^-\lambda(de)$.
\end{prop}

\noindent\emph{Proof.} As the parameters of the reflected BSDE \eqref{ekv:rbsde-jmp} satisfy the conditions for the comparison result of \cite{Quenez2014}, everything but the viscosity solution property follows by results presented therein. Existence of a unique viscosity solution in $\Pi^g_c$ to \eqref{ekv:obst-prob-n} can be shown as a special case of the method described in \cite{qvi-rbsde} (see, in particular, Theorem 4.5 therein). We can now argue as in \cite{HamMor16} and let $v\in\Pi^g_c$ be defined as $v(t,x):=Y^{t,x}_t$. Then, it can be shown (see Proposition 3.1 in \cite{HamMor16}) that $V^{t,x}_s(e)=v(s,X^{t,x}_s+\gamma(s,X^{t,x}_s,e))-v(s,X^{t,x}_s)$, $d\Prob\otimes ds\otimes \lambda(de)$-a.e. Moreover, the BSDE
\begin{align}\label{ekv:rbsde-n-2}
  \begin{cases}
     \tilde Y^{t,x}_s=\psi(X^{t,x}_T)+\int_s^T \tilde f^n(r,X^{t,x}_r, \tilde Y^{t,x}_r, \tilde Z^{t,x}_r, v(r,X^{t,x}_r+\gamma(r,X^{t,x}_r,\cdot))-v(r,X^{t,x}_r))dr
    \\
    \quad -\int_s^T  \tilde Z^{t,x}_r dW_r - \int_{s}^T\!\!\!\int_E  \tilde V^{t,x}_{r}(e)\mu(dr,de)+(\tilde K^{+,t,x}_T- \tilde K^{+,t,x}_s),\quad \forall s\in [t,T],
    \\
    \tilde Y^{t,x}_s\geq h(s,X^{t,x}_s),\, \forall s\in [t,T] \quad\text{and}\quad\int_t^T\!\! \big(\tilde Y^{t,x}_s-h(s,X^{t,x}_s)\big)dK^{+,t,x}_s,
  \end{cases}
\end{align}
admits a unique solution and $\tilde v(t,x):=\tilde Y^{t,x}_t$ belongs to $\Pi^g_c$ and is the unique viscosity solution to
\begin{align*}
\begin{cases}
  \min\{\tilde v(t,x)-h(t,x),-\tilde v_t(t,x)-\mcL \tilde v(t,x)+\mcK^n v(t,x)\\
  \quad-\tilde f(t,x,\tilde v(t,x),\sigma^\top(t,x)\nabla_x \tilde v(t,x))\}=0,\quad\forall (t,x)\in[0,T)\times \R^d \\
  \tilde v(T,x)=\psi(x).
\end{cases}
\end{align*}
On the other hand, $(Y^{t,x},Z^{t,x},V^{t,x},K^{+,t,x})$ is the unique solution to \eqref{ekv:rbsde-n-2} and we conclude that $\tilde v =v $ is the unique solution to \eqref{ekv:obst-prob-n}.\qed\\

To emphasize the dependence of the solution to \eqref{ekv:rbsde-jmp} on the parameter $n$, we will henceforth employ the notation $(Y^{t,x,n},Z^{t,x,n},V^{t,x,n},K^{+,t,x,n})$ for these quadruples.

\subsection{Preliminary estimates}
For $\nu\in\mcV_t$ we let $\E^\nu$ be expectation with respect to the probability measure $\Prob^\nu$ on $(\Omega,\mcF)$ defined by $d\Prob^\nu:=\kappa^\nu_T d\Prob$ with
\begin{align*}
\kappa^{\nu}_s&:=\mcE_s\Big(\int_{t}^\cdot\!\!\int_E(\nu_r(e)-1)(\mu(dr,de)-\lambda(de)dr)\Big)
\\
&:=\exp\Big(\int_{t}^s\!\!\!\int_E(1-\nu_r(e))\lambda(de)dr\Big)\prod_{t<\eta_j\leq s}\nu_{\eta_j}(\beta_j),
\end{align*}
where the sequence $(\eta_j,\beta_j)_{j\geq 1}$ is the one that appears in the Dirac decomposition $\mu=\sum_{j\geq 1}\delta_{(\eta_j,\beta_j)}$.

\begin{lem}\label{lem:SDEmoment}
Under Assumption~\ref{ass:onSDE}, the SDE \eqref{ekv:SDE-w-jmp} admits a unique solution for each $(t,x)\in[0,T]\times\R^d$. Furthermore, the solution has moments of all orders, in particular, for each $p\geq 0$, there is a constant $C>0$ such that
\begin{align}\label{ekv:SDEmoment}
\E^\nu\Big[\sup_{s\in [\zeta,T]}|X^{t,x}_s|^{p}\Big|\mcF^t_\zeta\Big]\leq C(1+|X^{t,x}_\zeta|^p),
\end{align}
$\Prob$-a.s.~for all $(t,x)\in [0,T]\times\R^n$, $\nu\in\mcV_t$ and $\zeta\in [t,T]$.
\end{lem}

\noindent\emph{Proof.} The proof is rather elementary and follows a similar structure to the proof of Proposition 4.2 in \cite{Perninge2022} (albeit the latter is confined to a Brownian filtration). It is provided in its entirety since some of its intermediate results are utilized later on. For each $j\in\bbN$, we let $X^{j}$ be the unique solution to the SDE
\begin{align*}
X^{j}_s=x + \int_t^s a(r,X^{j}_r)dr + \int_t^s \sigma(r,X^{j}_r)dW_r+\int_{t}^s\!\!\!\int_E  \ett_{[\mu((0,r),E)<j]}\gamma(s,X^{j}_{s-},e)\mu(dr,de),\quad\forall s\in [t,T],
\end{align*}
and note that, since $\mu$ has finite intensity, \rev{(outside of a $\Prob$-null set) we have $X^j\to X^{t,x}$ pointwisely as $j\to\infty$}. By Assumption~\ref{ass:onSDE}.(\ref{ass:onSDE-Gamma}) we get for $s\in [\eta_{j},T]$, using integration by parts, that
\begin{align*}
|X^{j}_s|^2&= |X^{j}_{\zeta\vee\eta_{j}}|^2+2\int_{(\zeta\vee\eta_{j})+}^s X^{j}_{r}dX^{j}_r+\int_{(\zeta\vee\eta_{j})+}^s d[X^{j},X^{j}]_r
\\
&\leq K^2_\Gamma\vee |X^{{j-1}}_{\zeta\vee\eta_{j}}|^2+2\int_{(\zeta\vee\eta_{j})+}^s X^{j}_{r} dX^{j}_r+\int_{(\zeta\vee\eta_{j})+}^s d[X^{j},X^{j}]_r.
\end{align*}
Now, either $|X^{{j-1}}_{\eta_{j}}|\leq K_\Gamma$ in which case
\begin{align*}
|X^{j}_s|^2\leq |X^{j}_{\zeta}|^2\vee K^2_\Gamma+2\int_{(\zeta\vee\eta_{j})+}^s X^{j}_{r} dX^{j}_r+\int_{(\zeta\vee\eta_{j})+}^s d[X^{j},X^{j}]_r.
\end{align*}
or $|X^{{j-1}}_{\eta_{j}}|> K_\Gamma$ implying that
\begin{align*}
|X^{j}_s|^2&\leq K^2_\Gamma\vee |X^{{j-2}}_{\zeta\vee\eta_{j-1}}|^2+2\int_{(\zeta\vee\eta_{j-1})+}^{\eta_{j}} X^{j-1}_{r} dX^{j-1}_r+\int_{(\zeta\vee\eta_{j-1})+}^{\eta_{j}} d[X^{j-1},X^{j-1}]_r
\\
&\quad+2\int_{(\zeta\vee\eta_{j})+}^s X^{j}_{r} dX^{j}_r+\int_{(\zeta\vee\eta_{j})+}^s d[X^{j},X^{j}]_r.
\end{align*}
In the latter case the same argument can be repeated and we conclude that
\begin{align}\label{ekv:X2-bound}
|X^{j}_s|^2&\leq |X^{j}_{\zeta}|^2\vee K_\Gamma^2+\sum_{i=j_0}^{j} \Big\{2\int_{(\zeta\vee\tilde\eta_{i})+}^{s\wedge\tilde\eta_{i+1}} X^{i}_{r}dX^{i}_r+\int_{(\zeta\vee\tilde\eta_{i})+}^{s\wedge\tilde\eta_{i+1}} d[X^{i},X^{i}]_r\Big\},
\end{align}
where $\tilde\eta_0=-1$, $\tilde\eta_i=\eta_i$ for $i=1,\ldots,j$ and $\tilde\eta_{j+1}=\infty$ and $j_0:=\max\{i\in \{1,\ldots,j\}:|X^{{i-1}}_{\eta_{i}}|\leq K_\Gamma\}\vee 0$. \rev{Here, $j_0$ is $\mcF_T$-measurable. This will, however, not render any issues as we take the supremum below.}

Now, since $X^{i}$ and $X^{j}$ coincide on $[0,\eta_{i+1\wedge j+1})$ we have
\begin{align*}
\sum_{i=j_0}^{j}\int_{(\zeta\vee\tilde\eta_{i})+}^{s\wedge\tilde\eta_{i+1}} X^{i}_{r} dX^{i}_r
&=\int_{\zeta\vee\eta_{j_0}}^s X^{j}_{r}a(r,X^{j}_r)dr+\int_{\zeta\vee\eta_{j_0}}^{s}X^{j}_{r}\sigma(r,X^{j}_r)dW_r,
\end{align*}
and
\begin{align*}
\sum_{i=j_0}^{j} \int_{(\zeta\vee\tilde\eta_{i})+}^{s\wedge\tilde\eta_{i+1}} d[X^{i},X^{i}]_r&=\int_{\zeta\vee\eta_{j_0}}^{s} \sigma^2(r,X^{j}_r)dr.
\end{align*}
Inserted in \eqref{ekv:X2-bound} this gives that
\begin{align}\nonumber
|X^{j}_s|^2&\leq |X^{j}_{\zeta}|^2\vee K_\Gamma^2+\int_{\eta_{j_0}}^s (2X^{j}_{s}a(r,X^{j}_r)+\sigma^2(r,X^{j}_r))dr+2\int_{\eta_{j_0}}^{s}X^{j}_{r}\sigma(r,X^{j}_r)dW_r
\\
&\leq |X^{j}_{\zeta}|^2+C\Big(1+\int_{\zeta}^{s}|X^{j}_{r}|^2dr + \sup_{\eta\in[\zeta,s]}\Big|\int_{\zeta}^{\eta}X^{j}_r\sigma(r,X^{j}_r)dW_r\Big|\Big)\label{ekv:X-bound1}
\end{align}
for all $s\in [\zeta,T]$. The Burkholder-Davis-Gundy inequality and the fact that the right-hand side of \eqref{ekv:X-bound1} does not depend on $\mu$ now gives that for any $\nu\in\mcV_t$ and $p\geq 2$,
\begin{align*}
\E^\nu\Big[\sup_{r\in[\zeta,s]}|X^{j}_r|^{p}\Big|\mcF^t_\zeta \Big]\leq |X^{j}_{\zeta}|^2+C\big(1+\E^\nu\Big[\int_{\zeta}^{s}|X^{j}_{r}|^{p}dr+\big(\int_{\zeta}^{s}|X^{j}_r|^4 dr\big)^{p/4}\Big|\mcF^t_\zeta\Big]\big).
\end{align*}
We can thus apply Gr\"onwall's lemma to conclude that for $p\geq 4$,
\begin{align*}
\E^\nu\Big[\sup_{s\in[\zeta,T]}|X^{j}_s|^{p}\Big|\mcF^t_\zeta\Big]&\leq C(1+ |X^{j}_{\zeta}|^{p}),
\end{align*}
$\Prob$-a.s., where the constant $C=C(T,p)$ does not depend on $\nu$ and $j$ and \eqref{ekv:SDEmoment} follows by letting $j\to\infty$ on both sides and using Fatou's lemma. The result for general $p\geq 0$ is then a simple consequence of Jensen's inequality.\qed\\

\begin{lem}\label{lem:Yn-bnd}
There is a $C>0$ such that
\begin{align}\label{ekv:vfs-bound}
|Y^{t,x,n}_s|\leq C(1+|X^{t,x}_s|^q)
\end{align}
for all $n\in\bbN$, $(t,x)\in[0,T]\times\R^d$ and $s\in[t,T]$.
\end{lem}

\noindent\emph{Proof.} We have
\begin{align}
h(s,X^{t,x}_s)\leq Y^{t,x,n}_s \leq Y^{t,x,0}_s
\end{align}
and the assertion follows by the polynomial growth assumptions on $h$ and the fact that $v_0\in\Pi^g_c$ (see Proposition~\ref{prop:rbsde-jmp}).\qed\\

%

We also make use of the following lemma which is given without proof as it follows immediately from the definitions:
\begin{lem}\label{lem:mcM-monotone}
Let $u,v:[0,T]\times \R^n\to\R$ be locally bounded functions. $\mcM$ is monotone (if $u\leq v$ pointwise, then $\mcM u\leq \mcM v$). Moreover, $\mcM(u_*)$ (resp. $\mcM(u^*)$) is l.s.c.~(resp. u.s.c.).
\end{lem}
In particular, it follows that $\mcM v$ is jointly continuous whenever $v$ is.\\

\section{The local setting\label{sec:local}}
Our approach to obtain existence of a unique viscosity solution to \eqref{ekv:var-ineq} goes through a fixed point argument. Before we are ready to proceed with this argument we need to solve the problem in the case where the driver $f$ is a local function. In this regard, we consider the PDE
\begin{align}\label{ekv:var-ineq-simp}
\begin{cases}
  \min\{v(t,x)-h(t,x),\max\{v(t,x)-\mcM v(t,x),-v_t(t,x)-\mcL v(t,x)\\
  \quad-\tilde f(t,x,v(t,x),\sigma^\top(t,x)\nabla_x v(t,x))\}\}=0,\quad\forall (t,x)\in[0,T)\times \R^d \\
  v(T,x)=\psi(x),
\end{cases}
\end{align}
where $\tilde f$ satisfies Assumption~\ref{ass:tilde-coeff}.

The arguments utilized in this section are based on penalization and use the unique viscosity solution, $v_n$, to \eqref{ekv:obst-prob-n}. By a comparison result for solutions to RBSDEs with jumps (see \eg Theorem 4.1 in \cite{Quenez2014}), we find that $(v_n)_{n\in\bbN}$ is a non-increasing sequence in $\Pi^g_c$ and by Lemma \ref{lem:Yn-bnd} this sequence is bounded from below in the sense that there is a constant $C>0$ such that $v_n(t,x)\geq -C(1+|x|^\rho)$ for all $(t,x)\in[0,T]\times\R^d$ and $n\in\bbN$. It then immediately follows that there is an upper semi-continuous function $v\in\Pi^g$ such that $v_n\searrow v$, pointwisely. We prove that $v$ is the unique viscosity solution to \eqref{ekv:var-ineq-simp} within the set of functions of polynomial growth.

We first show that $v$ satisfies the requirements for a viscosity solution at time $T$.

\begin{lem}\label{lem:simp-terminal}
For each $x\in\R^d$ it holds that $v^*(T,x)\leq \psi(x)$ and $v_*(T,x)\geq \psi(x)$.
\end{lem}

\noindent\emph{Proof.} First note that
\begin{align*}
  v^*(T,x)&\leq \limsup_{(t',x')\to(T,x) }v_{0}(t',x')=\psi(x)
\end{align*}
by continuity of $v_0$, proving the first inequality.

We turn to the second one which requires more work. In search for a contradiction, we assume that there is an $x'_0\in\R^d$ such that $v_*(T,x'_0) < \psi(x'_0)$. The proof is based on the following observation:
\begin{enumerate}[a)]
  \item If $v_*(T,x'_0) < \psi(x'_0)$ for some $x'_0\in\R^d$, then there is a (possible different) point $x_0\in\R^d$ such that $v_*(T,x_0) < \psi(x_0)\wedge \mcM v_*(T,x_0)$.
\end{enumerate}
To see this, we note that if \emph{a)} does not hold, then by lower semi-continuity, there is an $e'_0\in E$ such that
\begin{align*}
  v_*(T,x'_0)\geq \mcM v_*(T,x'_0)=v_*(T,x'_1)+\chi(T,x'_0,e'_0),
\end{align*}
with $x'_1:=x'_0+\gamma(T,x'_0,e'_0)$. Now, by Assumption~\ref{ass:on-coeff}.(\ref{ass:on-coeff-@end}), $\psi(x'_0) \leq \psi(x'_1) + \chi(T,x'_0,e'_0)$ and we find that
\begin{align*}
  v_*(T,x'_1) - \psi(x'_1)&\leq v_*(T,x'_0)-\chi(T,x'_0,e'_0) - \psi(x'_1)
  \\
  &\leq v_*(T,x'_0)-\chi(T,x'_0,e'_0) - (\psi(x'_0) - \chi(T,x'_0,e'_0))
  \\
  &< 0
\end{align*}
Hence, if \emph{a)} does not hold, then we must have $v_*(T,x'_1)\geq \mcM v_*(T,x'_1)$. We can repeat this argument indefinitely to find that, if \emph{a)} does not hold, then there is a sequence $(x'_j,e'_j)_{j=0}^\infty$ in $\R^d\times E$ such that $x'_{j+1}=x'_{j}+\gamma(T,x'_{j},e'_{j})$ and $v_*(T,x'_{j})\geq v_*(T,x'_{j+1})+\chi(T,x'_j,e'_j)$ for all $j\in\bbN$. \rev{We conclude that for each $0\leq k\leq l$, it holds that
\begin{align*}
  v_*(T,x'_k) &\geq  v_*(T,x'_l)+\sum_{j=k}^{l-1}\chi(T,x'_j,e'_{j}).
\end{align*}
Now, $(x'_j)_{j\geq 0}$ is a bounded sequence and thus contains a subsequence that converges to a point $\hat x\in\R^d$. Let $k_0:=\inf\{j\in\bbN:|x'_j-\hat x|\leq \delta_1(T)/2\}$ and let $k_i:=\inf\{j>k_{i-1}:|x'_j-\hat x|\leq \delta_1(T)/2\}$ for $i\geq 1$, then $|x'_{k_{i}}-x_{k_{i-1}}|\leq \delta_1(T)$ and the no-free-$\eps$-loop condition gives that
\begin{align*}
  v_*(T,x'_{k_0}) &\geq  v_*(T,x'_{k_i})+\sum_{j=k_0}^{k_{i}-1}\chi(T,x'_j,e'_j).
  \\
  &\geq h(T,x'_{k_i}) +i\delta_2(T)
\end{align*}
for each $i\in\bbN$,} a contradiction since $h$ is uniformly bounded from below on $[0,T]\times \Lambda_f(|x'_0|)$ and $v_*(T,x'_{k_0})\leq v_0(T,x'_{k_0})$.

We thus assume that there is a point $x_0$ and an $\eps>0$ such that
\begin{align}\label{ekv:contradiction}
  v_*(T,x_0)\leq \psi(x_0) \wedge \mcM v_*(T,x_0)-3\eps.
\end{align}
There is a sequence $(t_j,x_j,n_j)_{j\geq 0}$ in $[0,T)\times\R^d\times \bbN$ such that $(t_j,x_j)\to (T,x_0)$ and $v_{n_j}(t_j,x_j)\to v_*(T,x_0)$. In particular, if \eqref{ekv:contradiction} holds, then there is a $j_0$ such that
\begin{align*}
  v_{n_j}(t_j,x_j)\leq \psi(x_0)\wedge \mcM v_*(T,x_0)-2\eps
\end{align*}
whenever $j\geq j_0$. On the other hand, continuity of $\psi$ and lower semi-continuity of $\mcM v_*$ then imply that there is a $\delta'>0$ such that
\begin{align*}
  v_{n_j}(t_j,x_j)&\leq \psi(x)\wedge \mcM v_*(t,x)-\eps
  \\
  &\leq \psi(x)\wedge \mcM v_{n_j}(t,x)-\eps
\end{align*}
whenever $j\geq j_0$ and $(t,x)\in B_{\delta'}(T,x_0)\cap ([0,T]\times\R^d)$ (where $B_{\delta'}(T,x_0)$ is the open ball in $\R^{d+1}$ of radius $\delta'$ centered at $(T,x_0)$).

We introduce the stopping times
\begin{align*}
\theta_j:=\inf\{s\geq t_j:v_{n_j}(s,X^{t_j,x_j}_s)\geq \psi(X^{t_j,x_j}_s)\wedge \mcM v_{n_j}(s,X^{t_j,x_j}_s)\}
\end{align*}
and
\begin{align*}
\vartheta_j:=\inf\{s\geq t_j:(s,X^{t_j,x_j}_s)\notin B_{\delta'}(T,x_0)\text{ or }\mu((t_j,s],E)\geq 1\}.
\end{align*}
A standard dynamic programming result now gives that
\begin{align*}
     v_{n_j}(t_j,x_j)&= v_{n_j}(\theta_j\wedge \vartheta_j,X^{t_j,x_j}_{\theta_j\wedge \vartheta_j})+\int_{t_j}^{\theta_j\wedge \vartheta_j} \tilde f^{n_j}(r,X^{t_j,x_j}_r, Y^{t_j,x_j,n_j}_r,  Z^{t_j,x_j,n_j}_r, V^{t_j,x_j,n_j}_{r})dr
    \\
    &\quad-\int_{t_j}^{\theta_j\wedge \vartheta_j}   Z^{t_j,x_j,n_j}_r dW_r-\int_{t_j}^{\theta_j\wedge \vartheta_j}\!\!\!\int_E  V^{t_j,x_j,n_j}_{r}(e)\mu(dr,de)+K^{+,t_j,x_j,n_j}_{\theta_j\wedge \vartheta_j}- K^{+,t_j,x_j,n_j}_{t_j}.
\end{align*}
Since $V^{t_j,x_j,n_j}_{r}(e)=v_{n_j}(r,X^{t_j,x_j}_r+\gamma(r,X^{t_j,x_j}_r,e))-v_{n_j}(r,X^{t_j,x_j}_r)$, $d\Prob\otimes dr\otimes \lambda(de)$-a.e. (see Proposition 3.1 in \cite{HamMor16}) we get that,
\begin{align*}
    &\int_{t_j}^{\theta_j}\!\!\int_E\big(V^{t_j,x_j,n_j}_{r}(e)+\chi(r,X^{t_j,x_j}_r,e)\big)^-\lambda(de)dr
        \\
    &= \int_{t_j}^{\theta_j}\!\!\int_E\big(v_{n_j}(r,X^{t_j,x_j}_r +\gamma(r,X^{t_j,x_j}_r,e))+ \chi(r,X^{t_j,x_j}_r,e) - v_{n_j}(r,X^{t_j,x_j}_r)\big)^-\lambda(de)dr
    \\
    &\leq \int_{t_j}^{\theta_j}\!\!\int_E\big(\mcM v_{n_j}(r,X^{t_j,x_j}_r)-v_{n_j}(r,X^{t_j,x_j}_r)\big)^-\lambda(de)dr
    \\
    &=0.
\end{align*}
Consequently,
\begin{align*}
  \int_{t_j}^{\theta_j\wedge \vartheta_j} \tilde f^{n_j}(r,X^{t_j,x_j}_r, Y^{t_j,x_j,n_j}_r,  Z^{t_j,x_j,n_j}_r, V^{t_j,x_j,n_j}_{r})dr=\int_{t_j}^{\theta_j\wedge \vartheta_j} \tilde f(r,X^{t_j,x_j}_r, Y^{t_j,x_j,n_j}_r,  Z^{t_j,x_j,n_j}_r)dr
\end{align*}
and we conclude that
\begin{align*}
     v_{n_j}(t_j,x_j)&= v_{n_j}(\theta_j\wedge \vartheta_j,X^{t_j,x_j}_{\theta_j\wedge \vartheta_j})+\int_{t_j}^{\theta_j\wedge \vartheta_j} \tilde f(r,X^{t_j,x_j}_r, Y^{t_j,x_j,n_j}_r,  Z^{t_j,x_j,n_j}_r)dr
    \\
    &\quad-\int_{t_j}^{\theta_j\wedge \vartheta_j}   Z^{t_j,x_j,n_j}_r dW_r-\int_{t_j}^{\theta_j\wedge \vartheta_j}\!\!\!\int_E  V^{t_j,x_j,n_j}_{r}(e)\mu(dr,de)+K^{+,t_j,x_j,n_j}_{\theta_j\wedge \vartheta_j}- K^{+,t_j,x_j,n_j}_{t_j}.
\end{align*}
Standard arguments for reflected BSDEs now give that there is a $C>0$ such that
\begin{align*}
     \|Y^{t_j,x_j,n_j}\|_{\mcS^2_{[t,\theta_j\wedge \vartheta_j]}}+\|Z^{t_j,x_j,n_j}\|_{\mcH^2_{[t,\theta_j\wedge \vartheta_j]}(W)}+\|V^{t_j,x_j,n_j}\|_{\mcH^2_{[t,\theta_j\wedge \vartheta_j]}(\mu)}\leq C(1+|x_j|^\rho)
\end{align*}
for all $j\in\bbN$, and we find that
\begin{align*}
     v_{n_j}(t_j,x_j)&\geq \E\big[\ett_{[\theta_j<\vartheta_j]}v_{n_j}(\theta_j,X^{t_j,x_0}_{\theta_j})\big]-C\E\big[\ett_{[\theta_j\geq \vartheta_j]}\big]^{1/2}(1+|x_j|^\rho)-C(T-t_j)^{1/2}(1+|x_j|^\rho),
\end{align*}
where the constant $C>0$ can be chosen independent of $j$.  Now, since $v_{n_j}(T,x) = \psi(x)$, we have $\theta_j\leq T$, $\Prob$-a.s.~and on the subset of $\Omega$ where $[\theta_j<\vartheta_j]$ we have $v_{n_j}(\theta_j,X^{t_j,x_0}_{\theta_j})-v_{n_j}(t_j,x_j)\geq\eps$. We thus conclude that
\begin{align*}
    \eps\Prob[\theta_j< \vartheta_j]\leq C(T-t_j)^{1/2}(1+|x_j|^\rho)+C\E\big[\ett_{[\theta_j\geq \vartheta_j]}\big]^{1/2}(1+|x_j|^\rho),
\end{align*}
whenever $j\geq j_0$. On the other hand, $\Prob[\theta_j\geq \vartheta_j]\to 0$ as $j\to\infty$ and sending $j\to\infty$ gives the sought contradiction.\qed\\

\begin{thm}\label{thm:simp}
The function $v$ is continuous and thus belongs to $\Pi^g_c$. Moreover, it is the unique viscosity solution to \eqref{ekv:var-ineq-simp} in $\Pi^g$.
\end{thm}

\noindent\emph{Proof.} We prove that $v$ is a viscosity solution to \eqref{ekv:var-ineq}, then continuity and uniqueness will follow from the comparison principle stated in Proposition~\ref{app:prop:comp-visc} of Appendix \ref{app:uni-simp}. By continuity of $v_n$,
we have (see \eg \cite{Barles94}, p. 91),
\begin{align*}
  v_*(t,x)&=\liminf_{n\to\infty}\!{}_*v_n(t,x):=\liminf_{(n,t',x')\to(\infty,t,x) }v_n(t',x'),
  \\
   v(t,x)=v^*(t,x)&=\limsup_{n\to\infty}\!{}^*v_n(t,x):=\limsup_{(n,t',x')\to(\infty,t,x) }v_n(t',x').
\end{align*}
\emph{Subsolution-property.} We first show that $v=v^*$ is a viscosity subsolution. For $(t,x)\in [0,T)\times\R^d$ and $(p,q,M)\in\bar J^+v(t,x)$ there is (by Lemma 6.1 in~\cite{UsersGuide}) a sequence $(t_j,x_j)_{j\geq 0}$ in $[0,T)\times\R^d$ and sequences $n_j\to\infty$ and $(p_j,q_j,M_j)_{j\geq 0}$ with $(p_j,q_j,M_j)\in J^+v_{n_j}(t_j,x_j)$ such that
\begin{align*}
  (t_j,x_j,v_{n_j}(t_j,x_j),p_j,q_j,M_j)\to (t,x,v(t,x),p,q,M).
\end{align*}
As $(p_j,q_j,M_j)\in J^+v_{n_j}(t_j,x_j)$ and $v_{n_j}$ is a viscosity subsolution to \eqref{ekv:obst-prob-n} with $n=n_j$, it holds that
\begin{align*}
  &\min\{v_{n_j}(t_j,x_j)-h(t_j,x_j),-p_j-<a(t_j,x_j),q_j>-\frac{1}{2}\trace(\sigma\sigma^\top(t_j,x_j)M_j)
  \\& + n_j\int_E(v_{n_j}(t_j,x_j+\gamma(t_j,x_j,e))+\chi(t_j,x_j,e)-v_{n_j}(t_j,x_j))^-\lambda(de)
  \\
  &-\tilde f(t_j,x_j,v_{n_j}(t_j,x_j),\sigma^\top(t_j,x_j)q_j)\}\leq 0.
\end{align*}
Now, whenever $v(t,x)>h(t,x)$ there is a $j_0\in\bbN$, such that $v_{n_j}(t_j,x_j)>h(t_j,x_j)$ for all $j\geq j_0$. Hence,
\begin{align}\nonumber
  &-p_j-<a(t_j,x_j),q_j>-\frac{1}{2}\trace(\sigma\sigma^\top(t_j,x_j)M_j)
  \\& +n_j\int_E(v_{n_j}(t_j,x_j+\gamma(t_j,x_j,e))+\chi(t_j,x_j,e)-v_{n_j}(t_j,x_j))^-\lambda(de)\nonumber
  \\
  &-\tilde f(t_j,x_j,v_{n_j}(t_j,x_j),\sigma^\top(t_j,x_j)q_j)\leq 0\label{ekv:subsolu}
\end{align}
whenever $j\geq j_0$. Sending $j\to\infty$ gives that
\begin{align*}
  -p\,-<a(t,x),q>-\frac{1}{2}\trace(\sigma\sigma^\top(t,x)M)-\tilde f(t,x,v(t,x),\sigma^\top(t,x)q)\leq 0.
\end{align*}
Assume that $v(t,x)>\inf_{e\in E}\{v(t,x+\gamma(t,x,e))+\chi(t,x,e)\}$. Then there exists an $e_0\in E$ such that
\begin{align*}
  v(t,x+\gamma(t,x,e_0))+\chi(t,x,e_0) - v(t,x)<0
\end{align*}
and by upper semi-continuity of $v$ there is an $\eps>0$ and an open neighborhood $E_0\in\mcB(E)$ of $e_0$ in $E$ such that
\begin{align*}
  v(t,x+\gamma(t,x,e))+\chi(t,x,e) - v(t,x)\leq -2\eps
\end{align*}
for all $e\in \bar E_0$. On the other hand, there is a sequence $(e_j)_{j\geq 0}$ such that
\begin{align*}
  \sup_{e\in \bar E_0}\{v_{n_j}(t_j,x_j+\gamma(t_j,x_j,e))+\chi(t_j,x_j,e)\} = v_{n_j}(t_j,x_j+\gamma(t_j,x_j,e_j))+\chi(t_j,x_j,e_j)
\end{align*}
for each $j\in\bbN$ and we can extract a subsequence $(\tilde n_j,\tilde t_j,\tilde x_j,\tilde e_j)_{j\geq 0}$ of $(n_j,t_j,x_j,e_j)_{j\geq 0}$ such that $\tilde e_j\to\tilde e\in \bar E_0$ and
\begin{align*}
  \limsup_{j\to\infty}\sup_{e\in \bar E_0}\{v_{n_j}(t_j,x_j+\gamma(t_j,x_j,e))+\chi(t_j,x_j,e)\} & = \lim_{j\to\infty} v_{\tilde n_j}(\tilde t_j,\tilde x_j+\gamma(\tilde t_j,\tilde x_j,\tilde e_j))+\chi(\tilde t_j,\tilde x_j,\tilde e_j)
  \\
  &
  \leq \lim_{j\to\infty} v_{n'}(\tilde t_j,\tilde x_j+\gamma(\tilde t_j,\tilde x_j,\tilde e_j))+\chi(\tilde t_j,\tilde x_j,\tilde e_j)
  \\
  &= v_{n'}(t,x+\gamma(t,x,\tilde e))+\chi(t,x,\tilde e)
\end{align*}
for any $n'\in\bbN$. Since the right-hand side converges to
\begin{align*}
  v(t,x+\gamma(t,x,\tilde e))+\chi(t,x,\tilde e)\leq \sup_{e\in \bar E_0}\{v(t,x+\gamma(t,x,e))+\chi(t,x,e)\}\leq v(t,x) -2\eps
\end{align*}
as $n'\to\infty$, we conclude that
\begin{align*}
  v_{n_j}(t_j,x_j+\gamma(t,x_j,e))+\chi(t_j,x_j,e) - v_{n_j}(t_j,x_j)\leq -\eps
\end{align*}
for all $e\in E_0$ and all $j$ sufficiently large. Consequently,
\begin{align*}
  \int_E(v_{n_j}(t_j,x_j+\gamma(t_j,x_j,e))+\chi(t_j,x_j,e)-v_{n_j}(t_j,x_j))^-\lambda(de)\geq \eps\lambda(E_0)
\end{align*}
for $j$ sufficiently large. However, since $\lambda$ has full topological support, $\lambda(E_0)>0$ and this contradicts the fact that \eqref{ekv:subsolu} holds for all $j$. We thereby conclude that $v(t,x)\leq \inf_{e\in E}\{v(t,x+\gamma(t,x,e))+\chi(t,x,e)\}$.

\bigskip

\emph{Supersolution-property.}  We turn to the supersolution property of $v_*$. For $(t,x)\in [0,T)\times\R^d$ and $(p,q,M)\in\bar J^-v_*(t,x)$ there is (again by Lemma 6.1 in~\cite{UsersGuide}) a sequence $(t_j,x_j)_{j\geq 0}$ in $[0,T)\times\R^d$ and sequences $n_j\to\infty$ and $(p_j,q_j,M_j)_{j\geq 0}$ with $(p_j,q_j,M_j)\in J^-v_{n_j}(t_j,x_j)$ such that
\begin{align*}
  (t_j,x_j,v_{n_j}(t_j,x_j),p_j,q_j,M_j)\to (t,x,v_*(t,x),p,q,M).
\end{align*}
As $(p_j,q_j,M_j)\in J^-v_{n_j}(t_j,x_j)$ and $v_{n_j}$ is a viscosity supersolution to \eqref{ekv:obst-prob-n} with $n=n_j$, it holds that
\begin{align}\nonumber
  &\min\{v_{n_j}(t_j,x_j)-h(t_j,x_j),-p_j-<a(t_j,x_j),q_j>-\frac{1}{2}\trace(\sigma\sigma^\top(t_j,x_j)M_j)
  \\& + n_j\int_E(v_{n_j}(t_j,x_j+\gamma(t_j,x_j,e))+\chi(t_j,x_j,e)-v_{n_j}(t_j,x_j))^-\lambda(de)\nonumber
  \\
  &-\tilde f(t_j,x_j,v_{n_j}(t_j,x_j),\sigma^\top(t_j,x_j)q_j)\}\geq 0\label{ekv:subsolu-j}
\end{align}
In particular, $v_{n_j}(t_j,x_j)\geq h(t_j,x_j)$ and it follows by continuity of $h$ that $v_*(t,x)\geq h(t,x)$. Suppose now that $v_*(t,x)<\inf_{e\in E}\{v_*(t,x+\gamma(t,x,e))+\chi(t,x,e)\}$ implying the existence of an $\eps>0$ such that
\begin{align*}
v_*(t,x)\leq v_*(t,x+\gamma(t,x,e))+\chi(t,x,e)-\eps,\quad\forall e\in E.
\end{align*}
Since $v_{n_j}(t_j,x_j)\to v_{*}(t,x)$ and
\begin{align*}
&\liminf_{j\to\infty}\inf_{e\in\E}\{v_{n_j}(t_j,x_j+\gamma(t_j,x_j,e)) + \chi(t_j,x_j,e)\}
\\
&\geq \liminf_{j\to\infty}\inf_{e\in\E}\{v_{*}(t_j,x_j+\gamma(t_j,x_j,e)) + \chi(t_j,x_j,e)\}
\\
&\geq \inf_{e\in\E}\{v_{*}(t,x+\gamma(t,x,e)) + \chi(t,x,e)\},
\end{align*}
by lower semi-continuity of $\mcM v_*$ (see Lemma~\ref{lem:mcM-monotone}), this implies that
\begin{align*}
v_{n_j}(t,x)\leq v_{n_j}(t_j,x_j+\gamma(t_j,x_j,e))+\chi(t_j,x_j,e),\quad\forall e\in E,
\end{align*}
whenever $j$ is sufficiently large. In particular, we get that
\begin{align*}
  \int_E(v_{n_j}(t_j,x_j+\gamma(t_j,x_j,e))+\chi(t_j,x_j,e)-v_{n_j}(t_j,x_j))^-\lambda(de)=0
\end{align*}
for $j$ sufficiently large. Taking the limit in \eqref{ekv:subsolu-j} we thus find that in this situation
\begin{align*}
  &-p\,-<a(t,x),q>-\frac{1}{2}\trace(\sigma\sigma^\top(t,x)M)-\tilde f(t,x,v_{*}(t,x),\sigma^\top(t,x)q)\geq 0,
\end{align*}
proving the supersolution property of $v_*$.\qed\\

Much like the functions $v_n$, for each $(t,x)\in [0,T]\times\R^d$, the processes $(Y^{t,x,n})_{n\in\bbN}$ form a non-increasing sequence that is bounded from below by $h(\cdot,X^{t,x}_\cdot)\in\mcS^2$, implying the existence of an $\bbF^t$-progressively measurable process $Y^{t,x}$ such that $Y^{t,x,n}\searrow Y^{t,x}$, pointwisely. Taking the limit on both sides of $v_n(\eta,X^{t,x}_\eta)=Y^{t,x,n}_\eta$ gives that $Y^{t,x}_\eta=v(\eta,X^{t,x}_\eta)$ for each $\eta\in\mcT_t$. In particular, $Y^{t,x}$ is \cadlag and thus belongs to $\mcS^2_t$. In our pursuit of a related optimal stopping problem, we let for each $\tau\in\mcT_t$, the process $Y^{t,x,\tau}$ be the first component in the quadruple of processes $(Y^{t,x,\tau},Z^{t,x,\tau},V^{t,x,\tau},K^{-,t,x,\tau}) \in \mcS^2_{[0,\tau]} \times\mcH^2_{[0,\tau]}(W) \times\mcH^2_{[0,\tau]}(\mu) \times \mcS^2_{[0,\tau],i}$ that is the unique maximal solution to the BSDE with constrained jumps,
\begin{align}\label{ekv:stopped-bsde-simp}
  \begin{cases}
     Y^{t,x,\tau}_s=\Psi(\tau,X^{t,x}_\tau)+\int_s^\tau \tilde f(r,X^{t,x}_r, Y^{t,x,\tau}_r, Z^{t,x,\tau}_r)dr -\int_s^\tau  Z^{t,x,\tau}_r dW_r-\int_{s}^\tau\!\!\!\int_E  V^{t,x,\tau}_{r}(e)\mu(dr,de)
    \\
    \quad-(K^{-,t,x,\tau}_\tau- K^{-,t,x,\tau}_s),\quad \forall s\in [0,\tau],
    \\
    V^{t,x,\tau}_s(e)\geq - \chi(s,X^{t,x}_{s-},e),\quad d\Prob\otimes ds\otimes \lambda(de)-a.e.
  \end{cases}
\end{align}
Theorem 3.1 of \cite{imp-stop-game} seamlessly transfers to the present context, allowing us to conclude the following:

\begin{prop}\label{prop:Y-rep-simp}
For each $(t,x)\in [0,T]\times\R^d$ and $\eta\in\mcT_t$, the process $Y^{t,x}$ can be represented as
\begin{align}\label{ekv:opt-stop-simp}
  Y^{t,x}_\eta=\esssup_{\tau\in\mcT_\eta}Y^{t,x,\tau}_\eta=\essinf_{\nu\in\mcV}\esssup_{\tau\in\mcT_\eta}P^{t,x;\tau,\nu}_\eta,
\end{align}
where given $\nu\in\mcV$ and $\tau\in\mcT$, the triple $(P^{t,x;\tau,\nu},Q^{t,x;\tau,\nu},S^{t,x;\tau,\nu})\in\mcS^2_{[0,\tau]}\times\mcH^2_{[0,\tau]}(W)\times\mcH^2_{[0,\tau]}(\mu)$ is the unique solution to the standard BSDE
\begin{align}\nonumber
  P^{t,x;\nu,\tau}_s&=\Psi(\tau,X^{t,x}_\tau)+\int_s^\tau \tilde f^\nu(r,X^{t,x}_r,P^{t,x;\nu,\tau}_r,Q^{t,x;\nu,\tau}_r,S^{t,x;\nu,\tau}_{r})dr
  \\
  &\quad-\int_s^\tau Q^{t,x;\nu,\tau}_r dW_r-\int_{s}^\tau\!\!\!\int_E  S^{t,x;\nu,\tau}_{r}(e)\mu(dr,de),\quad\forall s\in [0,\tau]\label{ekv:non-ref-bsde-simp}
\end{align}
with driver
\begin{align*}
\tilde f^\nu(t,x,y,z,v):=\tilde f(t,x,y,z)+\int_E(v(e)+\chi(t,x,e))\nu_t(e)\lambda(de).
\end{align*}
\end{prop}

\bigskip

\begin{rem}\label{rem:Y-rep-conv}
For our purpose in the next section, instead of \eqref{ekv:opt-stop-simp}, we will utilize the slightly less general representation
\begin{align}\label{ekv:opt-stop-simp-conv}
  Y^{t,x}_s=\esssup_{\tau\in\mcT^t_s}Y^{t,x,\tau}_s=\essinf_{\nu\in\mcV_t}\esssup_{\tau\in\mcT^t_s}P^{t,x;\tau,\nu}_s,
\end{align}
which will turn out to be more convenient.
\end{rem}

\section{The general setting\label{sec:general}}
We now turn to the general setting of a non-local driver. Existence will again follow by an approximation routine and for $(t,x)\in[0,T]\times\R^d$ and $k\in\bbN$, we let $Y^{t,x,k}_s:=\esssup_{\tau\in\mcT^t_s}Y^{t,x,k;\tau}_s$ for all $s\in[t,T]$, where the quadruple $(Y^{t,x,k;\tau},Z^{t,x,k;\tau},V^{t,x,k;\tau},K^{-,t,x,k;\tau}) \in \mcS^2_{[t,\tau]} \times\mcH^2_{[t,\tau]}(W) \times\mcH^2_{[t,\tau]}(\mu) \times \mcS^2_{[t,\tau],i}$ is the unique maximal solution to
\begin{align}\label{ekv:rbsde_k}
\begin{cases}
  Y^{t,x,k;\tau}_s=\Psi(\tau,X^{t,x}_{\tau})+\int_s^{\tau} f(r,X^{t,x}_r,\bar Y^{k-1}(r,\cdot),Z^{t,x,k;\tau}_r)dr-\int_s^{\tau} Z^{t,x,k;\tau}_rdW_r
  \\
  \quad-\int_{s}^\tau\!\!\!\int_E  V^{t,x,k;\tau}_{r}(e)\mu(dr,de) - (K^{-,t,x,k;\tau}_{\tau}-K^{-,t,x,k;\tau}_s),\quad\forall s\in[t,{\tau}]
  \\
  V^{t,x,k;\tau}_s(e)\geq-\chi(s,X^{t,x}_s,e),\quad d\Prob\otimes ds\otimes\lambda(de)-a.e.,
\end{cases}
\end{align}
for $k\geq 1$, with $\bar Y^{k-1}(t,x):=Y^{t,x,k-1}_t$ and $\bar Y^{0}\equiv 0$.

\begin{prop}
There is a sequence $(\bar Y^{k}\in \Pi^g_c)_{k\geq 0}$ that satisfies the recursion above.
\end{prop}

\emph{Proof.} We need to show that for each $k\geq 1$, there is a $v_{k-1}\in\Pi^g_c$ such that $Y^{t,x,k-1}_t=v_{k-1}(t,x)$ for all $(t,x)\in [0,T]\times\R^n$. However, for $k=1$ this is immediate by the definition. Now, the result follows by using Remark~\ref{prop:Y-rep-simp}, Theorem~\ref{thm:simp} and induction.\qed\\

For $(t,\zeta)\in[0,T]\times\R_+$ and $\alpha\in\mcA_t$, we let $(\Upsilon^{t,\zeta;\alpha},\Theta^{t,\zeta;\alpha})\in\mcS^2_t\times\mcS^2_{t,i}$ solve the one-dimensional reflected SDE
\begin{align}\label{ekv:rsde}
\begin{cases}
  \Upsilon^{t,\zeta;\alpha}_s=\zeta^2\vee K_\Gamma^2+(4C_{a,\sigma}+2C_{a,\sigma}^2)\int_t^s(1+\Upsilon^{t,\zeta;\alpha}_r)dr+4C_{a,\sigma}\int_t^s (1+\Upsilon^{t,\zeta;\alpha}_r)\alpha_r dW_r+ \Theta^{t,\zeta;\alpha}_s,
  \\
  \quad\forall s\in[t,T], \\
  \Upsilon^{t,\zeta;\alpha}_s\geq \zeta^2\vee K_\Gamma^2\text{ and } \int_t^T(\Upsilon^{t,\zeta;\alpha}_s-(\zeta^2\vee K_\Gamma^2))d\Theta^{t,\zeta;\alpha}_s=0.
\end{cases}
\end{align}
We then set $R^{t,\zeta;\alpha}_s:=\sqrt{\Upsilon^{t,\zeta;\alpha}_s}$ and note that classically, we have
\begin{align*}
\E\Big[\sup_{s\in [t,T]}|R^{t,\zeta;\alpha}_s|^p\Big]\leq C(1+|\zeta\vee K_\Gamma|^p),
\end{align*}
for all $p\geq 2$.

\begin{lem}\label{lem:R-bounds-X}
For each $(t,x)\in [0,T]\times \R^d$, there is an $\alpha\in\mcA_t$ such that $|X^{t,x}_s|\leq R^{t,|x|;\alpha}_s$ for all $s\in [t,T]$, $\Prob$-a.s.
\end{lem}

\noindent\emph{Proof.} Since
\begin{align*}
|2x a(r,x)+\sigma^2(r,x)|\leq (4C_{a,\sigma}+2C_{a,\sigma}^2)(1+|x|^2),
\end{align*}
it follows from \eqref{ekv:X-bound1} that we can always choose $\alpha\in\mcA_t$ such that
\begin{align*}
2C_{a,\sigma}(1+\Upsilon^{t,|x|,\xi;\alpha}_r)\alpha_r=X^{t,x}_{r}\sigma(r,X^{t,x}_r),\quad \forall r\in [t,T],
\end{align*}
and the statement holds by \eqref{ekv:X-bound1}.\qed\\

For each $\varphi\in\Pi^g_c$, we let $\bar Y^{\varphi}\in\Pi^g_c$ be defined as $\bar Y^{\varphi}(t,x)=\sup_{\tau\in\mcT^t}Y^{t,x,\varphi;\tau}_t$, where the quadruple $(Y^{t,x,\varphi;\tau},Z^{t,x,\varphi;\tau},V^{t,x,\varphi;\tau},K^{-,t,x,\varphi;\tau}) \in\mcS^2_{[t,\tau]}\times\mcH^2_{[t,\tau]}(W)\times\mcH^2_{[t,\tau]}(\mu)\times \mcS^2_{[t,\tau],i}$ is the unique maximal solution to
\begin{align}\label{ekv:rbsde-varphi}
\begin{cases}
  Y^{t,x,\varphi;\tau}_s=\Psi(\tau,X^{t,x}_{\tau})+\int_s^{\tau} f(r,X^{t,x}_r,\varphi(r,\cdot),Z^{t,x,\varphi;\tau}_r)dr-\int_s^{\tau} Z^{t,x,\varphi;\tau}_rdW_r
  \\
  \quad-\int_{s}^\tau\!\!\!\int_E  V^{t,x,\varphi;\tau}_{r}(e)\mu(dr,de) - (K^{-,t,x,\varphi;\tau}_{\tau}-K^{-,t,x,\varphi;\tau}_s),\quad\forall s\in[t,{\tau}]
  \\
  V^{t,x,k;\tau}_s(e)\geq-\chi(s,X^{t,x}_s,e),\,d\Prob\otimes ds\otimes\lambda(de)-a.e.,
\end{cases}
\end{align}
and note that letting $(P^{t,x,\varphi;\nu,\tau},Q^{t,x,\varphi;\nu,\tau},S^{t,x,\varphi;\nu,\tau}) \in\mcS^2_{[t,\tau]}\times\mcH^2_{[t,\tau]}(W)\times\mcH^2_{[t,\tau]}(\mu)$ solve
\begin{align}\nonumber
  P^{t,x,\varphi;\nu,\tau}_s&=\Psi(\tau,X^{t,x}_\tau)+\int_s^\tau f^\nu(r,X^{t,x}_r,P^{t,x,\varphi;\nu,\tau}_r,\varphi(r,\cdot),Q^{t,x,\varphi;\nu,\tau}_r,S^{t,x,\varphi;\nu,\tau}_{r})dr
  \\
  &\quad-\int_s^\tau Q^{t,x,\varphi;\nu,\tau}_r dW_r-\int_{s}^\tau\!\!\!\int_E  S^{t,x,\varphi;\nu,\tau}_{r}(e)\mu(dr,de),\label{ekv:non-ref-bsde}
\end{align}
where
\begin{align*}
f^\nu(t,x,g,z,v):=f(t,x,g,z,v)+\int_E(v(e)+\chi(t,x,e))\nu_t(e)\lambda(de),
\end{align*}
Remark~\ref{rem:Y-rep-conv} gives that
\begin{align}\label{ekv:bar-Y-rep}
  \bar Y^{\varphi}(t,x)=\inf_{\nu\in\mcV_t}\sup_{\tau\in\mcT^t}P^{t,x,\varphi;\nu,\tau}_t.
\end{align}
The above definitions allow us to present the following result, which extends prior results derived in \cite{qvi-rbsde} (see Proposition 4.3 therein) and forms the basis for our contraction argument:
\begin{prop}\label{prop:contraction}
There is a $\kappa>0$ such that for all $\varphi,\tilde \varphi\in \Pi^g_c$ and $\zeta>0$, we have
\begin{align}
  &\sup_{\alpha\in\mcA^W}\E\Big[\int_0^Te^{\kappa t}\sup_{x\in\Lambda_f(R^{0,\zeta;\alpha}_t)}|\bar Y^{\tilde \varphi}(t,x)-\bar Y^{\varphi}(t,x)|^2dt\Big]\nonumber
  \\
  &\leq \frac{1}{4}\sup_{\alpha\in\mcA^W}\E\Big[\int_0^T e^{\kappa t}\sup_{x\in\Lambda_f(R^{0,\zeta;\alpha}_t)}|\tilde \varphi(t,x)-\varphi(t,x)|^2dt\Big].\label{ekv:int-contraction}
\end{align}
Furthermore, there is a $C>0$ such that
\begin{align}\label{ekv:sup-contraction}
  \sup_{t\in[0,T]}\sup_{x\in\Lambda_f(\zeta)}|\bar Y^{\tilde \varphi}(t,x)-\bar Y^{\varphi}(t,x)|^2\leq C\sup_{\alpha\in\mcA^W}\E\Big[\int_0^T \sup_{x\in\Lambda_f(R^{0,\zeta;\alpha}_t)}|\tilde \varphi(t,x)-\varphi(t,x)|^2dt\Big]
\end{align}
for each $\zeta>0$.
\end{prop}

\noindent\emph{Proof.} By \eqref{ekv:bar-Y-rep} we find that
\begin{align*}
  \bar Y^{\varphi}(t,x)- \bar Y^{\tilde\varphi}(t,x)&\leq \sup_{\nu\in\mcV_t}\sup_{\tau\in\mcT^t}(P^{t,x,\varphi;\nu,\tau}_t-P^{t,x,\tilde \varphi;\nu,\tau}_t).
\end{align*}
Since a similar inequality holds in the opposite direction, we get that
\begin{align}\label{ekv:YbyP}
  |\bar Y^{\varphi}(t,x)- \bar Y^{\tilde\varphi}(t,x)|&\leq \sup_{\nu\in\mcV_t}\sup_{\tau\in\mcT^t} |P^{t,x,\varphi;\nu,\tau}_t-P^{t,x,\tilde \varphi;\nu,\tau}_t|.
\end{align}
For $(\nu,\tau)\in \mcV_t\times\mcT^t$, let $(P,Q,S):=(P^{t,x,\varphi;\nu,\tau},Q^{t,x,\varphi;\nu,\tau},S^{t,x,\varphi;\nu,\tau})$ and $(\tilde P,\tilde Q,\tilde S):=(P^{t,x,\tilde \varphi;\nu,\tau},Q^{t,x,\tilde \varphi;\nu,\tau},S^{t,x,\tilde \varphi;\nu,\tau})$ and note that for $\kappa>0$, It\^o's formula applied to $e^{\kappa \cdot}|\tilde P-P|^2$ gives
\begin{align*}
  &e^{\kappa t}|\tilde P_t-P_t|^2+\int_t^\tau e^{\kappa s}|\tilde Q_s-Q_s|^2ds+\int_t^\tau\!\!\int_E e^{\kappa s}|\tilde S_s(e)-S_s(e)|^2d\mu(ds,de)
  \\&=-2\int_t^\tau e^{\kappa s}(\tilde P_s-P_s)(\tilde Q_s-Q_s)dW_s - 2\int_t^\tau\!\!\int_E e^{\kappa s}(\tilde P_s-P_s)(\tilde S_s(e)-S_s(e))d\mu(ds,de)
  \\
  &\quad -\kappa \int_t^\tau e^{\kappa s}|\tilde P_s-P_s|^2ds + 2\int_t^\tau e^{\kappa s}(\tilde P_s-P_s)(f(s,X^{t,x}_s,\tilde \varphi(s,\cdot),\tilde Q_s)-f(s,X^{t,x}_s,\varphi(s,\cdot),Q_s))ds
  \\
  &\quad + 2\int_t^\tau\!\!\int_E e^{\kappa s}(\tilde P_s-P_s)(\tilde S(e)-S(e))\nu_t(e)\lambda(de).
\end{align*}
By assumption
\begin{align*}
  |f(s,X^{t,x}_s,\varphi(s,\cdot),Q_s)-f(s,X^{t,x}_s,\tilde \varphi(s,\cdot),\tilde Q_s)|\leq k_f(\sup_{x'\in \Lambda_f(|X^{t,x}_s|)}|\tilde \varphi(s,x')-\varphi(s,x')| + |\tilde Q_s-Q_s|).
\end{align*}
Hence, taking the expectation w.r.t.~the measure $\Prob^\nu$ and using inequalities $2Cxy\leq (Cx)^2+y^2$ and $2xy\leq x^2/\sqrt{\kappa}+\sqrt{\kappa}y^2$ gives
\begin{align*}
  e^{\kappa t}|\tilde P_t-P_t|^2&\leq (C^2+C\sqrt{\kappa}-\kappa)\E^\nu\Big[ \int_t^\tau e^{\kappa s}|\tilde P_s-P_s|^2ds \Big]
  \\
  &\quad+\frac{C}{\sqrt{\kappa}}\E^\nu\Big[\int_t^Te^{\kappa s}\sup_{x'\in \Lambda_f(|X^{t,x}_s|)}|\tilde \varphi(s,x')-\varphi(s,x')|^2ds \Big].
\end{align*}
Now, pick $\kappa_0>0$ such that $\kappa_0\geq C^2+C\sqrt{\kappa_0}$ and note that for each $\kappa\geq\kappa_0$, we have
\begin{align*}
  e^{\kappa t}|P_t-\tilde P_t|^2&\leq \frac{C}{\sqrt{\kappa}}\E^\nu\Big[\int_t^Te^{\kappa s}\sup_{x'\in \Lambda_f(|X^{t,x}_s|)}|\tilde \varphi(s,x')-\varphi(s,x')|^2ds \Big]
  \\
  &\leq \frac{C}{\sqrt{\kappa}}\sup_{\alpha\in\mcA_t}\E^\nu\Big[\int_t^Te^{\kappa s}\sup_{x'\in \Lambda_f(|R^{t,|x|;\alpha}_s|)}|\tilde \varphi(s,x')-\varphi(s,x')|^2ds \Big]
  \\
  &\leq \frac{C}{\sqrt{\kappa}}\sup_{\alpha\in\mcA^W_t}\E\Big[\int_t^Te^{\kappa s}\sup_{x'\in \Lambda_f(|R^{t,|x|;\alpha}_s|)}|\tilde \varphi(s,x')-\varphi(s,x')|^2ds \Big],
\end{align*}
where the first inequality follows from Lemma~\ref{lem:R-bounds-X} while the second one is due to the fact that the SDE in \eqref{ekv:rsde} does not depend on $\mu$ (see \eg Section 4.1 of \cite{Bandini18}). Since the right-hand side is non-decreasing in $|x|$ and independent of $\nu$ and $\tau$, \eqref{ekv:YbyP} now gives that
\begin{align}\label{ekv:contraction-intermediate}
  e^{\kappa t}\sup_{x\in\Lambda_f(\zeta)}|\bar Y^{\tilde \varphi}(t,x)-\bar Y^{\varphi}(t,x)|^2&\leq \frac{C}{\sqrt{\kappa}}\sup_{\alpha\in\mcA^W_t}\E\Big[\int_t^Te^{\kappa s}\sup_{x'\in \Lambda_f(R^{t,\zeta;\alpha}_s)}|\tilde \varphi(s,x')-\varphi(s,x')|^2ds \Big],
\end{align}
for any $\zeta\geq 0$. In particular, as both sides are continuous in $\zeta$ a standard dynamic programming argument gives that for any $\alpha_1\in \mcA^W$, we have
\begin{align*}
  \E\Big[e^{\kappa t}\sup_{x\in\Lambda_f(R^{0,\zeta;\alpha_1}_t)}|\bar Y^{\tilde \varphi}(t,x)-\bar Y^{\varphi}(t,x)|^2\Big]&\leq \frac{C}{\sqrt{\kappa}}\sup_{\alpha\in\mcA^W_t}\E\Big[\int_t^Te^{\kappa s}\sup_{x'\in \Lambda_f(R^{0,\zeta;\alpha_1\oplus_t\alpha}_s)}|\tilde \varphi(s,x')-\varphi(s,x')|^2ds \Big].
\end{align*}
Taking the supremum with respect to $\alpha_1$ on the right-hand side and once again relying on a standard dynamic programming argument gives that
\begin{align*}
  \E\Big[e^{\kappa t}\sup_{x\in\Lambda_f(R^{0,\zeta;\alpha_1}_t)}|\bar Y^{\tilde \varphi}(t,x)-\bar Y^{\varphi}(t,x)|^2\Big]&\leq \frac{C}{\sqrt{\kappa}}\sup_{\alpha\in\mcA^W}\E\Big[\int_0^Te^{\kappa s}\sup_{x'\in \Lambda_f(R^{0,\zeta;\alpha}_s)}|\tilde \varphi(s,x')-\varphi(s,x')|^2ds\Big].
\end{align*}
Integrating with respect to time and using Fubini's theorem, we find that
\begin{align*}
  \E\Big[\int_0^T e^{\kappa t}\sup_{x\in\Lambda_f(R^{0,\zeta;\alpha_1}_t)}|\bar Y^{\tilde \varphi}(t,x)-\bar Y^{\varphi}(t,x)|^2 dt\Big]&\leq \frac{CT}{\sqrt{\kappa}}\sup_{\alpha\in\mcA^W}\E\Big[\int_0^Te^{\kappa s}\sup_{x'\in \Lambda_f(R^{0,\zeta;\alpha}_s)}|\tilde \varphi(s,x')-\varphi(s,x')|^2ds\Big]
\end{align*}
after which taking the supremum with respect to $\alpha_1\in\mcA^W$ and choosing $\kappa\geq (4CT)^2\vee\kappa_0$ gives the first inequality. To get \eqref{ekv:sup-contraction} we note that comparison gives that $R^{0,\zeta;\alpha}_s\geq R^{t,\zeta;\alpha}_s$ for all $s\in [t,T]$ and $\alpha\in\mcA^W$. From \eqref{ekv:contraction-intermediate} we thus get that
\begin{align*}
  \sup_{x\in\Lambda_f(\zeta)}|\bar Y^{\tilde \varphi}(t,x)-\bar Y^{\varphi}(t,x)|^2&\leq C\sup_{\alpha\in\mcA^W}\E\Big[\int_0^T\sup_{x'\in \Lambda_f(R^{0,\zeta;\alpha}_s)}|\tilde \varphi(s,x')-\varphi(s,x')|^2ds \Big]
\end{align*}
from which \eqref{ekv:sup-contraction} is immediate since the right-hand side is independent of $t$.\qed\\

We now introduce the norm $\|\cdot\|_{\zeta}$ on the space of jointly continuous functions of polynomial growth, $\Pi^g_c$, defined as
\begin{align*}
  \|\varphi\|_{\zeta}:=\sup_{\alpha\in\mcA^W}\E\Big[\int_0^Te^{\kappa t}\sup_{x\in \Lambda_f(R^{0,\zeta;\alpha}_t)}|\varphi(t,x)|^2dt\Big]^{1/2},
\end{align*}
with $\kappa>0$ as in Proposition~\ref{prop:contraction} and note that under $\|\cdot\|_{\zeta}$, the map $\Phi:\Pi^g_c\to\Pi^g_c$ that maps $\varphi$ to $\bar Y^{\varphi}$ is a contraction.

\begin{cor}\label{cor:unif-bound}
There are constants $C>0$ and $p\geq 0$ such that $|\bar Y^{k}(t,x)|\leq C(1+|x|^p)$ for all $(t,x)\in[0,T]\times\R^d$ and all $k\geq 0$.
\end{cor}

\noindent\emph{Proof.} First, we note that \eqref{ekv:int-contraction} and the triangle inequality imply that
\begin{align*}
  \|\bar Y^{k}\|_{\zeta}&\leq  \|\bar Y^{k}-\bar Y^{k-1}\|_{\zeta}+\|\bar Y^{k-1}\|_{\zeta}\leq \frac{1}{2}\|\bar Y^{k-1}-\bar Y^{k-2}\|_{\zeta}+\|\bar Y^{k-1}\|_{\zeta}\leq \frac{1}{2^{k-1}}\|\bar Y^{1}-\bar Y^{0}\|_{\zeta}+\|\bar Y^{k-1}\|_{\zeta}.
\end{align*}
However, as a similar scheme holds for $\|\bar Y^{k-1}\|_{\zeta}$ and since $\bar Y^{0}\equiv 0$ we conclude that
\begin{align*}
  \|\bar Y^{k}\|_{\zeta}&\leq  \sum_{j=1}^k\frac{1}{2^{j-1}}\|\bar Y^{1}\|_{\zeta}\leq 2\|\bar Y^{1}\|_{\zeta}.
\end{align*}
On the other hand, as $\bar Y^{1}\in\Pi^g_c$ there are constants $C>0$ and $p\geq 2$ such that $|\bar Y^{1}(t,x)|\leq C(1+|x|^p)$ and we conclude that
\begin{align*}
  \|\bar Y^{1}\|^2_{\zeta}&=\sup_{\alpha\in\mcA^W}\E\Big[\int_0^Te^{\kappa t}\sup_{x\in \Lambda_f(R^{0,\zeta;\alpha}_t)}|\bar Y^{1}(t,x)|^2dt\Big]
  \\
  &\leq C\Big(1+\sup_{\alpha\in\mcA^W}\E\Big[\sup_{t\in [0,T]}|R^{0,\zeta;\alpha}_t|^{2p}\Big]\Big)
  \\
  &\leq C(1+|\zeta|^{2p})
\end{align*}
implying the existence of a $C>0$ such that $\|\bar Y^{k}\|_{\zeta}\leq C(1+|\zeta|^{p})$ for all $k\geq 0$. Now, \eqref{ekv:sup-contraction} gives that
\begin{align*}
  \sup_{t\in[0,T]}\sup_{x\in\Lambda_f(\zeta)}|\bar Y^{k}(t,x)-\bar Y^{1}(t,x)|^2&\leq C\sup_{\alpha\in\mcA^W}\E\Big[\int_0^T \sup_{x\in\Lambda_f(R^{0,\zeta;\alpha}_t)}|\bar Y^{k-1}(t,x)|^2dt\Big]
  \\
  &\leq C(1+|\zeta|^{2p})
\end{align*}
where the constants $C>0$ and $p\geq 2$ do not depend on $k$ and the desired bound follows.\qed\\

Letting $v_k(t,x):=\bar Y^{k}(t,x)$, Proposition~\ref{prop:contraction} and Corollary~\ref{cor:unif-bound} \rev{imply} that there is a $v\in \Pi^g$ such that for each $\zeta>0$ we have $\|v_k-v\|_{\zeta}\to 0$ as $k\to\infty$.

\begin{thm}\label{thm:non-local}
$v$ is the unique viscosity solution in $\Pi^g_c$ to \eqref{ekv:var-ineq}.
\end{thm}

\noindent \emph{Proof.} First, \eqref{ekv:sup-contraction} implies that the convergence is uniform on compact subsets of $[0,T]\times\R^n$ and since $v_k$ is jointly continuous for each $k\geq 0$, we conclude that $v$ is also jointly continuous. This in turn gives that $\Phi(v)$ is well defined and we conclude that $\Phi(v)=v$, establishing existence of a solution to the optimal stopping problem \eqref{ekv:opt-stop}-\eqref{ekv:SDE-w-jmp}. Moreover, if $\tilde Y$ is another solution, then $\tilde v(t,x):=\tilde Y^{t,x}_t$ must also satisfy $\Phi(\tilde v)=\tilde v$. However, then repeated use of the contraction property in \eqref{ekv:int-contraction} gives that $\|\tilde v-v\|_{\zeta}=0$ and by continuity we conclude that $\tilde v=v$ implying by uniqueness of the non-linear Snell envelope in \eqref{ekv:opt-stop-simp-conv}, as obtained in Proposition~\ref{prop:Y-rep-simp}, that $\tilde Y^{t,x}=Y^{t,x}$. The solution to the optimal stopping problem \eqref{ekv:opt-stop}-\eqref{ekv:SDE-w-jmp} is thus unique.

Utilizing, once more, the connection between the optimal stopping problem \eqref{ekv:opt-stop}-\eqref{ekv:SDE-w-jmp} and PDEs with obstacles we conclude that $v$ is a viscosity solution to \eqref{ekv:var-ineq}. Suppose now that there exists another function $\tilde v\in\Pi^g_c$ that solves \eqref{ekv:var-ineq} and let $\bar v=\Phi(\tilde v)$, then by Theorem~\ref{thm:simp} and Proposition~\ref{prop:Y-rep-simp} we conclude that $\bar v\in\Pi^g_c$ is the unique solution to
\begin{align*}
\begin{cases}
  \min\{\bar v(t,x)-h(t,x),\max\{\bar v(t,x)- \mcM \bar v(t,x),-\bar v_t(t,x)-\mcL \bar v(t,x)\\
  \quad-f(t,x,\tilde v(t,\cdot),\sigma^\top(t,x)\nabla_x \bar v(t,x))\}\}=0,\quad\forall (t,x)\in[0,T)\times \R^d \\
  \bar v(T,x)=\psi(x),
\end{cases}
\end{align*}
But then $\bar v=\tilde v$ implying that $\tilde v$ is a fixed point of $\Phi$ and since $v$ is the only fixed point of $\Phi$ in the set of jointly continuous functions of polynomial growth we must have $\tilde v=v$.\qed\\

\begin{cor}\label{cor:non-local}
$Y^v$ is the unique solution to the optimal stopping problem \eqref{ekv:opt-stop}-\eqref{ekv:SDE-w-jmp}.
\end{cor}

\appendix
\section{Uniqueness of viscosity solutions in the local framework\label{app:uni-simp}}
In this section we prove the critical comparison principle for the PDE with local driver,
\begin{align}\label{ekv:var-ineq-simp-A}
\begin{cases}
  \min\{v(t,x) - h(t,x),\max\{v(t,x)-\mcM v(t,x),-v_t(t,x)-\mcL v(t,x)\\
  \quad-\tilde f(t,x,v(t,x),\sigma^\top(t,x)\nabla_x v(t,x))\}\}=0,\quad\forall (t,x)\in[0,T)\times \R^d \\
  v(T,x)=\psi(x),
\end{cases}
\end{align}
that was treated in Section~\ref{sec:local} (see \eqref{ekv:var-ineq-simp}). We need the following lemma:

\begin{lem}\label{lem:is-super}
Let $v$ be a supersolution to \eqref{ekv:var-ineq-simp-A} satisfying
\begin{align*}
\forall (t,x)\in[0,T]\times\R^d,\quad |v(t,x)|\leq C(1+|x|^{2\varrho})
\end{align*}
for some $\varrho>0$. Then there is a $\varpi_0 > 0$ such that for any $\varpi>\varpi_0$ and $\theta > 0$, the function $v + \theta e^{-\varpi t}(1+((|x|-K_\Gamma)^+)^{2\varrho+2})$ is also a supersolution to \eqref{ekv:var-ineq-simp-A}.
\end{lem}

\noindent \emph{Proof.} With $w(t,x):=v(t,x) + \theta e^{-\varpi t}(1+((|x|-K_\Gamma)^+)^{2\varrho+2})$ we note that, since $v$ is a supersolution and $\theta e^{-\varpi T}(1+((|x|-K_\Gamma)^+)^{2\varrho+2})\geq 0$, we have
\begin{align*}
  w(t,x)\geq v(t,x)\geq \ett_{[t<T]}h(t,x)+\ett_{[t=T]}\psi(x)
\end{align*}
so that $w(t,x)\geq h(t,x)$ for all $(t,x)\in [0,T)\times\R^d$ and the terminal condition holds. Assume now that
\begin{align*}
  v(t,x) - \inf_{e\in E}\{v(t,x+\gamma(t,x,e))+\chi(t,x,e)\}\geq 0.
\end{align*}
In this case,
\begin{align*}
&w(t,x) - \inf_{e\in E}\{w(t,x+\gamma(t,x,e))+\chi(t,x)\}
\\
&=v(t,x) + \theta e^{-\varpi t}(1+((|x|-K_\Gamma)^+)^{2\varrho+2})
\\
&\quad - \inf_{e\in E}\{v(t,x+\gamma(t,x,e)) + \theta e^{-\varpi t}(1+((|x+\gamma(t,x,e)|-K_\Gamma)^+)^{2\varrho+2}) +\chi(t,x,e)\}
\\
&\geq v(t,x) - \inf_{e\in E}\{v(t,x+\gamma(t,x,e))+\chi(t,x,e)\}
\\
&\quad+\theta e^{-\varpi t}\{1+((|x|-K_\Gamma)^+)^{2\varrho+2}- \sup_{e\in E}(1+((|x+\gamma(t,x,e)|-K_\Gamma)^+)^{2\varrho+2})\}.
\end{align*}
Now, either $|x|\leq K_\Gamma$ in which case it follows by \eqref{ekv:imp-bound} that $|x+\gamma(t,x,e)|\leq K_\Gamma$ or $|x|> K_\Gamma$ and \eqref{ekv:imp-bound} gives that $|x+\gamma(t,x,e)|\leq |x|$. We conclude that
\begin{align*}
  w(t,x) - \inf_{e\in E}\{w(t,x+\gamma(t,x,e))+\chi(t,x,e)\}\geq 0.
\end{align*}
Consider instead the case when
\begin{align*}
  v(t,x) - \inf_{e\in E}\{v(t,x+\gamma(t,x,e))+\chi(t,x,e)\}< 0
\end{align*}
and let $\varphi\in C^{1,2}([0,T]\times\R^d\to\R)$ be such that $\varphi-w$ has a local maximum of 0 at $(t,x)$ with $t<T$. Then $(\tilde t,\tilde x)\mapsto\tilde \varphi(\tilde t,\tilde x):=\varphi (\tilde t,\tilde x)-\theta e^{-\varpi \tilde t}(1+((|\tilde x|-K_\Gamma)^+)^{2\varrho+2})\in C^{1,2}([0,T]\times\R^d\to\R)$ and $\tilde \varphi-v$ has a local maximum of 0 at $(t,x)$. Since $v$ is a viscosity supersolution, we have
\begin{align*}
    &-\partial_t(\varphi(t,x)-\theta e^{-\varpi t}(1+((|x|-K_\Gamma)^+)^{2\varrho+2}))-\mcL (\varphi(t,x)-\theta e^{-\varpi t}(1+((|x|-K_\Gamma)^+)^{2\varrho+2}))
    \\
    &-\tilde f(t,x,\varphi(t,x)-\theta e^{-\varpi t}(1+((|x|-K_\Gamma)^+)^{2\varrho+2}),\sigma^\top(t,x)\nabla_x (\varphi(t,x)-\theta e^{-\varpi t}(1+((|x|-K_\Gamma)^+)^{2\varrho+2})))\geq 0.
\end{align*}
Consequently,
\begin{align*}
&-\partial_t\varphi(t,x)-\mcL \varphi(t,x)-\tilde f(t,x,\varphi(t,x),\sigma^\top(t,x)\nabla_x \varphi(t,x))
\\
&\geq \theta \varpi e^{-\varpi t}(1+((|x|-K_\Gamma)^+)^{2\varrho+2})-\theta\mcL e^{-\varpi t}(1+((|x|-K_\Gamma)^+)^{2\varrho+2})
\\
&\quad \tilde f(t,x,\varphi(t,x)-\theta e^{-\varpi t}(1+((|x|-K_\Gamma)^+)^{2\varrho+2}),\sigma^\top(t,x)\nabla_x (\varphi(t,x)-\theta e^{-\varpi t}(1+((|x|-K_\Gamma)^+)^{2\varrho+2})))
\\
&\quad-\tilde f(t,x,\varphi(t,x),\sigma^\top(t,x)\nabla_x \varphi(t,x))
\\
&\geq \theta \varpi e^{-\varpi t}(1+((|x|-K_\Gamma)^+)^{2\varrho+2})-\theta C(1+\varrho) e^{-\varpi t}(1+((|x|-K_\Gamma)^+)^{2\varrho+2})
\\
&\quad -k_f(1+\varrho)\theta e^{-\varpi t}(1+((|x|-K_\Gamma)^+)^{2\varrho+2}),
\end{align*}
where the right hand side is non-negative for all $\theta> 0$ and all $\varpi>\varpi_0$ for some $\varpi_0>0$.\qed\\

We have the following result, the proof of which we omit since it is classical:
\begin{lem}\label{lem:integ-factor}
For any $\kappa\in\R$, a locally bounded function $v:[0,T]\times \R^d\to\R$ is a viscosity supersolution (resp. subsolution) to \eqref{ekv:var-ineq-simp-A} if and only if $\check v(t,x):=e^{\kappa t}v(t,x)$ is a viscosity supersolution (resp. subsolution) to
\begin{align}\label{ekv:var-ineq-if}
\begin{cases}
  \min\{\check v(t,x)-e^{\kappa t}h(t,x),\max\{\check v(t,x)-\inf_{e\in E}\{\check v(t,x+\gamma(t,x,e))+e^{\kappa t}\chi(t,x,e)\},-\check v_t(t,x)
  \\
+\kappa \check v(t,x)-\mcL \check v(t,x)-e^{\kappa t}\tilde f(t,x,e^{-\kappa t}\check v(t,x),e^{-\kappa t}\sigma^\top(t,x)\nabla_x \check v(t,x))\}=0,\quad\forall (t,x)\in[0,T)\times \R^d \\
  \check v(T,x)=e^{\kappa T}\psi(x).
\end{cases}
\end{align}
\end{lem}
\begin{rem}
Here, it is important to note that $\check h(t,x):=e^{\kappa t}h(t,x)$, $\check \chi(t,x,e):=e^{\kappa t}\chi(t,x,e)$, $\check f(t,x,y,z):=-\kappa y+e^{\kappa t}\tilde f(t,x,e^{-\kappa t}y,e^{-\kappa t}z)$ and $\check \psi(x):=e^{\kappa T}\psi(x)$ satisfy Assumption~\ref{ass:on-coeff}. In particular, this implies that Lemma~\ref{lem:is-super} holds for supersolutions to \eqref{ekv:var-ineq-if} as well.
\end{rem}

We have the following comparison result for viscosity solutions in $\Pi^g$:

\begin{prop}\label{app:prop:comp-visc}
Let $v$ (resp. $u$) be a supersolution (resp. subsolution) to \eqref{ekv:var-ineq-simp-A}. If $u,v\in \Pi^g$, then $u\leq v$.
\end{prop}

\noindent \emph{Proof.} First, we note that it is sufficient to show that the statement holds for solutions to~\eqref{ekv:var-ineq-if} for some $\kappa\in \R$. We thus assume that $v$ (resp.~$u$) is a viscosity supersolution (resp.~subsolution) to \eqref{ekv:var-ineq-if} for $\kappa\in \R$ specified below. Furthermore, we may without loss of generality assume that $v$ is l.s.c.~and $u$ is u.s.c.

By assumption, $u,v\in \Pi^g$, which implies that there are constants $C>0$ and $\varrho>0$ such that
\begin{align}\label{ekv:uv-bound}
|v(t,x)|+|u(t,x)|\leq C(1+|x|^{2\varrho}).
\end{align}
Now, for any $\varpi>0$ we only need to show that
\begin{align*}
w(t,x)&=w^{\theta,\varpi}(t,x):=v(t,x)+\theta e^{-\varpi t}(1+((|x|-K_\Gamma)^+)^{2\varrho+2})
\\
&\geq u(t,x)
\end{align*}
for all $(t,x)\in[0,T]\times\R^d$ and any $\theta>0$. Then the result follows by taking the limit $\theta\to 0$. We know from Lemma~\ref{lem:is-super} that there is a $\varpi_0>0$ such that $w$ is a supersolution to \eqref{ekv:var-ineq-if} for each $\varpi\geq\varpi_0$ and $\theta>0$. We thus assume that $\varpi\geq\varpi_0$.

We search for a contradiction and assume that there is a $(t_0,x_0)\in [0,T]\times \R^d$ such that $u(t_0,x_0)>w(t_0,x_0)$. By \eqref{ekv:uv-bound}, there is for each $\theta>0$ a $R> K_\Gamma$ such that
\begin{align*}
w(t,x)\geq u(t,x),\quad\forall (t,x)\in[0,T]\times\R^d,\:|x|\geq R.
\end{align*}
Our assumption thus implies that there is a point $(\bar t,\bar x)\in[0,T)\times B_R$ (the open unit ball of radius $R$ centered at 0) such that
\begin{align*}
\max_{(t,x)\in[0,T]\times\R^d}(u(t,x)-w(t,x))&=\max_{(t,x)\in[0,T)\times B_R}(u(t,x)-w(t,x))
\\
&=u(\bar t,\bar x)-w(\bar t,\bar x)=\eta>0.
\end{align*}
We first show that there is at least one point $(t^*,x^*)\in[0,T)\times B_R$ such that
\begin{enumerate}[a)]
  \item $u(t^*,x^*)-w(t^*,x^*)= \eta$,
  \item $u(t^*,x^*)>\check h(t^*,x^*)$ and
  \item $w(t^*,x^*)<\inf_{e\in E}\{v(t^*,x^*+\gamma(t^*,x^*,e))+\check \chi(t^*,x^*,e)\}$.
\end{enumerate}
Assume first that $u(t,x)\leq \check h(t,x)$ for some $(t,x)\in[0,T]\times\R^d$. Since $w$ is a supersolution, we have $w(t,x)\geq \check h(t,x)$ and it follows that $u(t,x)-w(t,x)\leq 0$ contradicting that $u(t,x)-w(t,x)= \eta$. In particular, any point satisfying \emph{a)} must also satisfy \emph{b)}.

We proceed by assuming that $w(t,x)\geq\inf_{e\in E}\{w(t,x+\gamma(t,x,e))-\check \chi(t,x,e)\}$ for all $(t,x)\in A:=\{(s,y)\in[0,T]\times\R^d: u(s,y)-w(s,y)=\eta\}$. Indeed, as $w$ is l.s.c. and $\gamma$ is continuous, there is an $e_1\in E$ such that
\begin{align}\label{ekv:equal}
w(\bar t,\bar x)\geq\inf_{e\in E}\{w(\bar t,\bar x+\gamma(\bar t,\bar x,e))+\check \chi(\bar t,\bar x,e)\}=w(\bar t,\bar x+\gamma(\bar t,\bar x,e_0))+\check \chi(\bar t,\bar x,e_0).
\end{align}
Now, set $x_1=\bar x+\gamma(\bar t,\bar x,e_0)$ and note that since
 \begin{align*}
|x+\gamma(t,x,e)|< R,\quad \forall(t,x,e)\in [0,T]\times B_R\times E
\end{align*}
it follows that $x_1\in B_R$. Moreover, as $u$ is a subsolution with $u(\bar t,\bar x) > \check h(\bar t,\bar x)$ it satisfies
\begin{align*}
  u(\bar t,\bar x) - (u(\bar t,\bar x+\gamma(\bar t,\bar x,e_0)) + \check \chi(\bar t,\bar x,e_0))\leq 0,
\end{align*}
that is
\begin{align*}
  u(\bar t,x_1))\geq u(\bar t,\bar x)-\check \chi(t,\bar x,e_0)
\end{align*}
and we conclude from \eqref{ekv:equal} that
\begin{align*}
  u(\bar t,x_1)- w(\bar t,x_1)&\geq u(\bar t,\bar x)-\check \chi(\bar t,\bar x,e_0)-(w(\bar t,\bar x)-\check \chi(t,\bar x,e_0))
  \\
  &=u(\bar t,\bar x)-w(\bar t,\bar x)=\eta.
\end{align*}
Hence, $(\bar t,x_1)\in A$ and by our assumption it follows that there is an $e_1\in E$ such that
\begin{align*}
u(\bar t,x_1)=u(\bar t,x_1+\gamma(\bar t,x_1,e_1))+\check \chi(\bar t,x_1,e_1)
\end{align*}
and a corresponding $x_2:=x_1+\gamma(\bar t,x_1,e_1)\in B_R$. \rev{Now, this process can be repeated indefinitely to find a sequence $(e_j)_{j\geq 0}$ in $E$ such that for any $0\leq k\leq l$ we have
\begin{align*}
  w(\bar t,x_k)\geq w(\bar t,x_l)+\sum_{j=k}^{l-1}\check \chi(\bar t,x_{j},e_{j}),
\end{align*}}
with $x_0:=\bar x$ and $x_j:=x_{j-1}+\gamma(\bar t,x_{j-1},e_{j-1})$. We can now repeat the argument in the proof of Lemma~\ref{lem:simp-terminal} to get a contradiction as $w(t,x)\geq \check h(t,x)$, where the latter is bounded on $[0,T]\times \bar B_R$. We can thus find a $(t^*,x^*)\in [0,T)\times B_R$ such that \emph{a)}-\emph{c)} above holds.

Since $\tilde f$ is Lipschitz in $y$ and $z$ for $(t,x)\in[0,T]\times\bar B_R$, the remainder of the proof follows along the lines of the proof of Proposition 4.1 in~\cite{Morlais13} (see also Proposition 6.4 in \cite{PerningeJMAA23}) and is included only for the sake of completeness.

Next, we assume without loss of generality that $\varrho\geq 2$ and define
\begin{align*}
\Phi_n(t,x,y):=u(t,x)-w(t,y)-\varphi_n(t,x,y),
\end{align*}
where
\begin{align*}
  \varphi_n(t,x,y):=\frac{n}{2}|x-y|^{2\varrho}+|x-x^*|^2+|y-x^*|^2+(t-t^*)^2.
\end{align*}
Since $u$ is u.s.c.~and $w$ is l.s.c.~there is a triple $(t_n,x_n,y_n)\in[0,T]\times \bar B_R\times \bar B_R$ (with $\bar B_R$ the closure of $B_R$) such that
\begin{align*}
  \Phi_n(t_n,x_n,y_n)=\max_{(t,x,y)\in [0,T]\times \bar B_R\times \bar B_R}\Phi_n(t,x,y).
\end{align*}
Now, the inequality $2\Phi_n(t_n,x_n,y_n)\geq \Phi_n(t_n,x_n,x_n)+\Phi_n(t_n,y_n,y_n)$ gives
\begin{align*}
n|x_n-y_n|^{2\varrho}\leq u(t_n,x_n)-u(t_n,y_n)+w(t_n,x_n)-w(t_n,y_n).
\end{align*}
Consequently, $n|x_n-y_n|^{2\varrho}$ is bounded (since $u$ and $w$ are bounded on $[0,T]\times\bar B_R\times\bar B_R$) and $|x_n-y_n|\to 0$ as $n\to\infty$. We can, thus, extract subsequences $n_l$ such that $(t_{n_l},x_{n_l},y_{n_l})\to (\tilde t,\tilde x,\tilde x)$ as $l\to\infty$. Since
\begin{align*}
u(t^*,x^*)-w(t^*,x^*)\leq \Phi_n(t_n,x_n,y_n)\leq u(t_n,x_n)-w(t_n,y_n),
\end{align*}
it follows that
\begin{align*}
u(t^*,x^*)-w(t^*,x^*)&\leq \limsup_{l\to\infty} \{u(t_{n_l},x_{n_l})-w(t_{n_l},y_{n_l})\}
\\
&\leq u(\tilde t,\tilde x)-w(\tilde t,\tilde x)
\end{align*}
and as the right-hand side is dominated by $u(t^*,x^*)-w(t^*,x^*)$ we conclude that
\begin{align*}
  u(\tilde t,\tilde x)-w(\tilde t,\tilde x)=u(t^*,x^*)-w(t^*,x^*).
\end{align*}
In particular, this gives that $\lim_{l\to\infty}\Phi_n(t_{n_l},x_{n_l},y_{n_l})=u(\tilde t,\tilde x)-w(\tilde t,\tilde x)$ which implies that
\begin{align*}
  \limsup_{l\to\infty} n_l|x_{n_l}-y_{n_l}|^{2\varrho}= 0
\end{align*}
and
\begin{align*}
  (t_{n_l},x_{n_l},y_{n_l})\to (t^*,x^*,x^*).
\end{align*}
Since $u$ is u.s.c. and $w$ is l.s.c., we must also have
\begin{align*}
  (u(t_{n_l},x_{n_l}),w(t_{n_l},y_{n_l}))\to (u(t^*,x^*),w(t^*,x^*)).
\end{align*}
We can thus extract a subsequence $(\tilde n_l)_{l\geq 0}$ of $(n_l)_{l\geq 0}$ such that $t_{\tilde n_l}<T$, $|x_{\tilde n_l}|\vee |y_{\tilde n_l}|<R$ and
\begin{align*}
  u(t_{\tilde n_l},x_{\tilde n_l})-w(t_{\tilde n_l},y_{\tilde n_l})\geq \frac{\eta}{2}
\end{align*}
for all $l\in\bbN$. Moreover, since $\check h$ is continuous and $(t,y)\mapsto \inf_{e\in E}\{w(t,y+\gamma(t,y,e))+\check \chi(t,y,e)\}$ is l.s.c.~(see Lemma~\ref{lem:mcM-monotone}), there is an $l_0\geq 0$ such that
\begin{align*}
  u(t_{\tilde n_l},x_{\tilde n_l})-\check h(t_{\tilde n_l},x_{\tilde n_l})>0,
\end{align*}
and
\begin{align*}
  w(t_{\tilde n_l},y_{\tilde n_l})-\inf_{e\in E}\{w(t_{\tilde n_l},y_{\tilde n_l}+\gamma(t_{\tilde n_l},y_{\tilde n_l},e))-\check \chi(t_{\tilde n_l},y_{\tilde n_l},e)\}<0,
\end{align*}
for all $l\geq l_0$. To simplify notation we will, from now on, denote $(\tilde n_l)_{l\geq l_0}$ simply by $n$.\\

By Theorem 8.3 of~\cite{UsersGuide} there are $(p^u_n,q^u_n,M^u_n)\in \bar J^{+}u(t_n,x_n)$ and $(p^w_n,q^w_n,M^w_n)\in \bar J^{-}w(t_n,y_n)$, such that
\begin{align*}
\begin{cases}
  p^u_n-p^w_n=\partial_t\varphi_n(t_n,x_n,y_n)=2(t_n-t^*)
  \\
  q^u_n=D_x\varphi_n(t_n,x_n,y_n)=n\varrho(x-y)|x-y|^{2\varrho-2}+2(x-x^*)
  \\
  q^w_n=-D_y\varphi_n(t_n,x_n,y_n)=n\varrho(x-y)|x-y|^{2\varrho-2}+2(x-x^*)
\end{cases}
\end{align*}
and for every $\epsilon>0$,
\begin{align*}
  \left[\begin{array}{cc} M^n_x & 0 \\ 0 & -M^n_y\end{array}\right]\leq B(t_n,x_n,y_n)+\epsilon B^2(t_n,x_n,y_n),
\end{align*}
where $B(t_n,x_n,y_n):=D^2_{(x,y)}\varphi_n(t_n,x_n,y_n)$. Now, we have
\begin{align*}
  D^2_{(x,y)}\varphi_n(t,x,y)=\left[\begin{array}{cc} D_x^2\varphi_n(t,x,y) & D^2_{yx}\varphi_n(t,x,y) \\ D^2_{xy}\varphi_n(t,x,y) & D_y^2\varphi_n(t,x,y)\end{array}\right]
  = \left[\begin{array}{cc} n\xi(x,y)+2 I & -n\xi(x,y) \\ -n\xi(x,y) & n\xi(x,y) +2 I \end{array}\right]
\end{align*}
where $I$ is the identity-matrix of suitable dimension and
\begin{align*}
  \xi(x,y):=\varrho|x-y|^{2\varrho-4}\{|x-y|^2I+2(\varrho-1)(x-y)(x-y)^\top\}.
\end{align*}
In particular, since $x_n$ and $y_n$ are bounded, choosing $\epsilon:=\frac{1}{n}$ gives that
\begin{align}\label{ekv:mat-bound}
  \tilde B_n:=B(t_n,x_n,y_n)+\epsilon B^2(t_n,x_n,y_n)\leq Cn|x_n-y_n|^{2\varrho-2}\left[\begin{array}{cc} I & -I \\ -I & I \end{array}\right]+C I.
\end{align}
By the definition of viscosity supersolutions and subsolutions we have that
\begin{align*}
&-p^u_n+\kappa u(t_n,x_n)-a^\top(t_n,x_n)q^u_n-\frac{1}{2}\trace [\sigma^\top(t_n,x_n)M^u_n\sigma(t_n,x_n)]
\\
&-e^{\kappa t_n}\tilde f(t_n,x_n,e^{-\kappa t_n}u(t_n,x_n),e^{-\kappa t_n}\sigma^\top(t_n,x_n)q^u_n)\leq 0
\end{align*}
and
\begin{align*}
&-p^w_n+\kappa w(t_n,y_n)-a^\top(t_n,y_n)q^w_n-\frac{1}{2}\trace [\sigma^\top(t_n,y_n)M^w_n\sigma(t_n,y_n)]
\\
&-e^{\kappa t_n}\tilde f(t_n,y_n,e^{-\kappa t_n}w(t_n,y_n),e^{-\kappa t_n}\sigma^\top(t_n,x_n)q^w_n)\geq 0.
\end{align*}
Combined, this gives that
\begin{align*}
\kappa (u(t_n,x_n)-w(t_n,y_n))&\leq p^u_n+a^\top(t_n,x_n)q^u_n+\frac{1}{2}\trace [\sigma^\top(t_n,x_n)M^u_n\sigma(t_n,x_n)]
\\
&+e^{\kappa t_n}\tilde f(t_n,x_n,e^{-\kappa t_n}u(t_n,x_n),e^{-\kappa t_n}\sigma^\top(t_n,x_n)q^u_n)
\\
&-p^w_n-a^\top(t_n,y_n)q^w_n-\frac{1}{2}\trace [\sigma^\top(t_n,y_n)M^w_n\sigma(t_n,y_n)]
\\
&-e^{\kappa t_n}\tilde f(t_n,y_n,e^{-\kappa t_n}w(t_n,y_n),e^{-\kappa t_n}\sigma^\top(t_n,x_n)q^w_n)
\end{align*}
Collecting terms we have that
\begin{align*}
p^u_n-p^w_n&=2(t_n-t^*)
\end{align*}
and since $a$ is Lipschitz continuous in $x$ and bounded on $\bar B_R$, we have
\begin{align*}
  a^\top(t_n,x_n)q^u_n-a^\top(t_n,y_n)q^w_n&\leq  (a^\top(t_n,x_n)-a^\top(t_n,y_n))n\varrho(x_n-y_n)|x_n-y_n|^{2\varrho-2}
  \\
  &\quad+C(|x_n-x^*|+|y_n-x^*|)
  \\
  &\leq C(n|x_n-y_n|^{2\varrho}+|x_n-x^*|+|y_n-x^*|),
\end{align*}
where the right-hand side tends to 0 as $n\to\infty$. Let $s_x$ denote the $i^{\rm th}$ column of $\sigma(t_n,x_n)$ and let $s_y$ denote the $i^{\rm th}$ column of $\sigma(t_n,y_n)$ then by the Lipschitz continuity of $\sigma$ and \eqref{ekv:mat-bound}, we have
\begin{align*}
  s_x^\top M^u_n s_x-s_y^\top M^w_n s_y&=\left[\begin{array}{cc} s_x^\top  & s_y^\top \end{array}\right]\left[\begin{array}{cc} M^u_n  & 0 \\ 0 &-M^w_n\end{array}\right]\left[\begin{array}{c} s_x \\ s_y \end{array}\right]
  \\
  &\leq \left[\begin{array}{cc} s_x^\top  & s_y^\top \end{array}\right]\tilde B_n\left[\begin{array}{c} s_x \\ s_y \end{array}\right]
  \\
  &\leq C(n|x_n-y_n|^{2\varrho}+|x_n-y_n|)
\end{align*}
and we conclude that
\begin{align*}
  \limsup_{n\to\infty}\frac{1}{2}\trace [\sigma^\top(t_n,x_n)M^u_n\sigma(t_n,x_n)-\sigma^\top(t_n,y_n)M^w_n\sigma(t_n,y_n)]\leq 0.
\end{align*}
Finally, we have that
\begin{align*}
  &e^{\kappa t_n}\tilde f(t_n,x_n,e^{-\kappa t_n}u(t_n,x_n),e^{-\kappa t_n}\sigma^\top(t_n,x_n)q^u_n)-e^{\kappa t_n}\tilde f(t_n,y_n,e^{-\kappa t_n}w(t_n,y_n),e^{-\kappa t_n}\sigma^\top(t_n,x_n)q^w_n)
  \\
  &\leq k_f(u(t_n,x_n)-w(t_n,y_n)+|\sigma^\top(t_n,x_n)q^u_n-\sigma^\top(t_n,x_n)q^w_n|)
  \\
  &\quad + e^{\kappa t_n}|\tilde f(t_n,x_n,e^{-\kappa t_n}u(t_n,x_n),e^{-\kappa t_n}\sigma^\top(t_n,x_n)q^u_n)-\tilde f(t_n,y_n,e^{-\kappa t_n}u(t_n,x_n),e^{-\kappa t_n}\sigma^\top(t_n,x_n)q^u_n)|
\end{align*}
Repeating the above argument and using that $\tilde f$ is jointly continuous in $(t,x)$ uniformly in $(y,z)$ we get that the upper limit of the right-hand side when $n\to\infty$ is bounded by $k_f(u(t_n,x_n)-w(t_n,y_n))$. Put together, this gives that
\begin{align*}
(\kappa-k_f) \limsup_{n\to\infty}(u(t_n,x_n)-w(t_n,y_n))&\leq 0
\end{align*}
and choosing $\kappa>k_f$ gives a contradition.\qed\\


\bibliographystyle{plain}
\bibliography{qvi-stop_ref}

\begin{thebibliography}{10}

\bibitem{Bandini18}
E.~Bandini, A.~Cosso, M.~Fuhrman, and H.~Pham.
\newblock Backward sdes for optimal control of partially observed
  path-dependent stochastic systems: a control randomization approach.
\newblock {\em Ann. Appl. Probab.}, 28(3):1634--1678, 2018.

\bibitem{Barles94}
G.~Barles.
\newblock {\em Solutions de viscosit\'e des \'equations de Hamilton-Jacobi},
  volume~17 of {\em Math\'ematiques et Applications}.
\newblock Springer, Paris, 1994.

\bibitem{Barles1997}
G.~Barles, R.~Buckdahn, and E.~Pardoux.
\newblock Backward stochastic differential equations and integral-partial
  differential equations.
\newblock {\em Stoch. Stoch. Rep.}, 60(1-2):57--83, 1997.

\bibitem{BensLionsImpulse}
A.~Bensoussan and J.L. Lions.
\newblock {\em Impulse Control and Quasivariational inequalities}.
\newblock Gauthier-Villars, Montrouge, France, 1984.

\bibitem{Bouchard09}
B.~Bouchard.
\newblock A stochastic target formulation for optimal switching problems in
  finite horizon.
\newblock {\em Stochastics}, 81:171--197, 2009.

\bibitem{UsersGuide}
M.~G. Crandall, H.~Ishii, and P.~L. Lions.
\newblock Users guide to viscosity solutions of second order partial
  differential equations.
\newblock {\em Bulletin of the American Mathematical Society}, 27(1):1--67,
  1992.

\bibitem{Dumitrescu15}
R.~Dumitrescu, M.-C. Quenez, and A.~Sulem.
\newblock Optimal stopping for dynamic risk measures with jumps and obstacle
  problems.
\newblock {\em J Optim Theory Appl}, 167:219--242, 2015.

\bibitem{Fuhrman2020}
M.~Fuhrman and M.~Morlais.
\newblock Optimal switching problems with an infinite set of modes: An approach
  by randomization and constrained backward sdes.
\newblock {\em Stochastic Process. Appl.}, 130:5(5):3120--3153, 2020.

\bibitem{Fuhrman15}
M.~Fuhrman and H.~Pham.
\newblock Randomized and backward sde representation for optimal control of
  non-markovian sdes.
\newblock {\em Ann. Appl. Probab.}, 25(4):2134--2167, 2015.

\bibitem{HamMnif2023}
S.~Hamad{\`e}ne, M.~Mnif, and S.~Neffati.
\newblock Viscosity solution of system of integro-partial differential
  equations with interconnected obstacles of non-local type without
  monotonicity conditions.
\newblock {\em J Dyn Diff Equat}, 35(2):1151--1173, 2023.

\bibitem{Morlais13}
S.~Hamad\`ene and M.~A. Morlais.
\newblock Viscosity solutions of systems of pdes with interconnected obstacles
  and switching problem.
\newblock {\em Appl Math Optim.}, 67:163--196, 2013.

\bibitem{HamMor16}
S.~Hamad\`ene and M.~A. Morlais.
\newblock Viscosity solutions for second order integro-differential equations
  without monotonicity condition: the probabilistic approach.
\newblock {\em Stochastics: An International Journal of Probability and
  Stochastic Processes}, 88(4):632--649, 2016.

\bibitem{KharroubiQVI}
I.~Kharroubi, J.~Ma, H.~Pham, and J.~Zhang.
\newblock Backward sdes with constrained jumps and quasi-variational
  inequalities.
\newblock {\em Ann. Probab.}, 38(2):794--840, 2010.

\bibitem{OksenSulemBok}
B.~{\O}ksendal and A.~Sulem.
\newblock {\em Applied Stochastic Control of Jump Diffusions}.
\newblock Springer, 2007.

\bibitem{Perninge2022}
M.~Perninge.
\newblock Sequential systems of reflected backward stochastic differential
  equations with application to impulse control.
\newblock {\em Appl Math Optim}, 86(19), 2022.

\bibitem{qvi-rbsde}
M.~Perninge.
\newblock Probabilistic representation of viscosity solutions to
  quasi-variational inequalities with non-local drivers.
\newblock {\em ESAIM Control Optim.\,Calc.\,Var.}, 25:1--28, 2023.

\bibitem{PerningeJMAA23}
M.~Perninge.
\newblock Zero-sum stochastic differential games of impulse versus continuous
  control by fbsdes.
\newblock {\em J. Math. Anal. Appl.}, 527, 2023.

\bibitem{imp-stop-game}
M.~Perninge.
\newblock Optimal stopping of bsdes with constrained jumps and related zero-sum
  games.
\newblock {\em Stochastic Process. Appl.}, 173, 2024.

\bibitem{Quenez2014}
M-C. Quenez and A.~Sulem.
\newblock Reflected bsdes and robust optimal stopping for dynamic risk measures
  with jumps.
\newblock {\em Stochastic Process. Appl.}, 124:3031--3054, 2014.

\end{thebibliography}
\end{document}